# Algebraic Kaprekar's routine architecture (II)

**F. Nuez**

**Retired professor of Universidad Politécnica de Valencia (Spain)**

**Abstract**

In general terms, we establish algebraic relations that numbers must satisfy in order for their images to match after one or several transformations. Some groups associated with these relationships are identified, such as the Klein group. Such equivalences are applied to numbers of 2, 3, 4 or 5 digits. The relationship between cycles and the transformation trees' structure are analyzed.

## 1. Introduction

Kaprekar's routine consists in sorting the digits of a number n in descending order, resulting X. Then, Y is obtained by sorting the digits in ascending order, and these numbers are subtracted n' = X-Y. When iterating the process with n' and beyond, Kaprekar (1949) showed that if n has four non-identical digits such iteration leads to 6174. Additionally, if the number n has three digits the routine leads to 495. These numbers are known as Kaprekar constants.

In the first part of this paper (Nuez, 2021) we have formalized the process by introducing parametric functions that allow to establish algebraic relationships along the transformation chains. This has allowed us to deduce constants and cycles from these functions. As application examples, the sets of numbers of 2, 3, 4, and 5 digits have been studied.

There are numerous works on the Kaprekar process (in references we show some of them). Recent studies are uncovering connections of this transformation with various sectors of Number Theory (Yamagami, 2018) and there are several authors who suspect there are important gaps in the understanding of the process.

In this paper we analyze the algebraic architecture of transformation trees. For this we develop a methodology based on binary equivalence relations that may be useful in the study of other transformation chains.

## 2. Notation

As a reminder, let us summarize some of the terminology defined and developed in Nuez (2021).

The reference sets $A_w$ are those consisting of natural numbers of up to w digits, excluding those with repeated digits. If the number of digits is lower than w, it is completed by adding zeros to its left. In this paper we study the cases w = 2, 3, 4 and 5. The absorbing subsets $B_w$ are defined as those consisting of all the images of $A_w$ numbers.

Each number has a parameter α associated to it, resulting from subtracting its extreme digits after sorting it in descending order, which is indicated as O(n). If the number has 4 or 5 digits a second parameter ß (ß≤α) is defined, resulting from the subtraction of those digits adjacent to those in O(n). These parameters characterize number n and determine what its image n' is. The expression p(n) = (α,ß) indicates the parameters of the number. There is the following fundamental relationship: p(m) = p(n) ↔ K(m) = K(n).

The transformation process is formalized by the following functions

- Functions f: p (n) = (α,ß),  f (α,ß) = n'
  They yield the numerical image of a number of especified parameters
  Thus, in $A_4$
  $f_1$ (α,ß) = (α  ß-1  9-ß  10-α), $0 < α \leq 9$, $0 < ß \leq 9$, $α \geq ß$
  $f_2$ (α,0) = (α-1  9  9  10-α), $0 < α \leq 9$
- Functions $K_i$:  p (n) = (α,ß), p (n') = (α',ß'),  $K_i$ (α,ß) = (α',ß')
  There are always several functions $K_i$, each of them acting on specific values of α and ß. All of these functions have been developed in the previous paper and are listed in the Annex.

In general notation, in order to refer to the "image of" symbol K is used. Its argument may be a number or the parameters of a number, resulting in numerical images or parametrical images, respectively. Thus,

K (n) = n';   K (α,ß) = (α',ß'),  p (n) = (α,ß), p (n') = (α',ß')

In the second case, we simply indicate that the transformation algorithm is operating, without specifying its domain of existence, unlike $K_i$. All numbers can be transformed with a unique image and, in consequence, so can their parameters.

The repeated application of operator K is indicated with a superscript,
K (0 0 0 1) = 0 9 9 9 , $K^2$ (0 0 0 1) = K [K(0 0 0 1)] = 8 9 9 1,  $K^5$ (0 0 0 1) = 6174
and for other numbers in Graph 1
K (2 5 0 9 ) = 9 2 6 1, $K^2$ (2 5 0 9) = 8 3 5 2, $K^3$ (2 5 0 9) = 6 1 7 4
K (3 0 9 5) = 9 1 7 1, $K^3$ (3 0 9 5) = 6174
To simplify we shall use $K^1$ (n) = K (n), and by convention, $K^0$ (n) = n.

## 3. Trees of parametric transformation

Since parameters unambiguously define the image, the graph consisting of the transformations of all numbers is considerably simplified if the parameters characterizing such numbers are represented. In the case of $A_4$ we go from 9990 numbers to 54 parametric combinations (examples in the rest of the paragraph refer to $A_4$).

Which function $K_i$ to use will depend on the values of α and ß, and it is indicated in Tables 2 and 3 of the previous paper (see the Annex). For example (9,3) satisfies the requirements of $K_2$: $1 < 3 \leq 5$, $11 \leq 9+3 \leq 14$, $1 \leq 9-3 \leq 7$, $6 \leq 9 \leq 9$ and its image will be $K_2$ (9,3) = (2x9-10, 10-2x3) = (8,4)

Some parameters α and ß can be transformed by several functions, all of which yield the same image. For example, (5,4) can be transformed by $K_5$, $K_6$, $K_{13}$, $K_{14}$ $K_{17}$ and $K_{18}$, whose image is (2,0).

The expression (α,ß) → (α’,ß’) or, equivalently, K (α,ß) = (α’,ß’), is an epimorphism of $A_4$ in $B_4$. Basic functions $f_1$ and $f_2$ determine that all numbers in $A_4$ whose parameters are α and ß turn into the same number of $B_4$, whose parameters are α’ and ß’. But different parameters can give the same parametric image. Therefore, in Graph 1 it is observed that 8991 and 1998 or other non-represented ones such as 8100 or 3084 all with parameters (8,1) yield the same image, 8082 of parameters (8,6). But there are other numbers, e.g. 9171 of parameters (8,6), which come from numbers with parameters other than (8,1). Is the case of 9171 of (9,2).

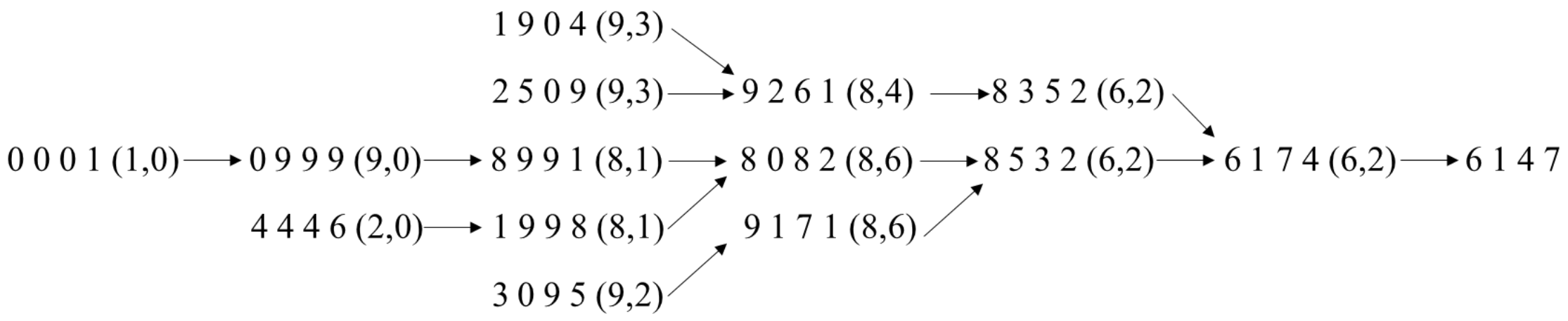

Graph 1. Example of numeric and parametric transformations

This Graph shows sequences of numbers which have not branched out in the direction of the transformation. It also contains parameters which characterize such sequences as a consequence of the image of any number being unique.

The whole tree of parametric transformations in $A_4$ is shown in Graph 2.

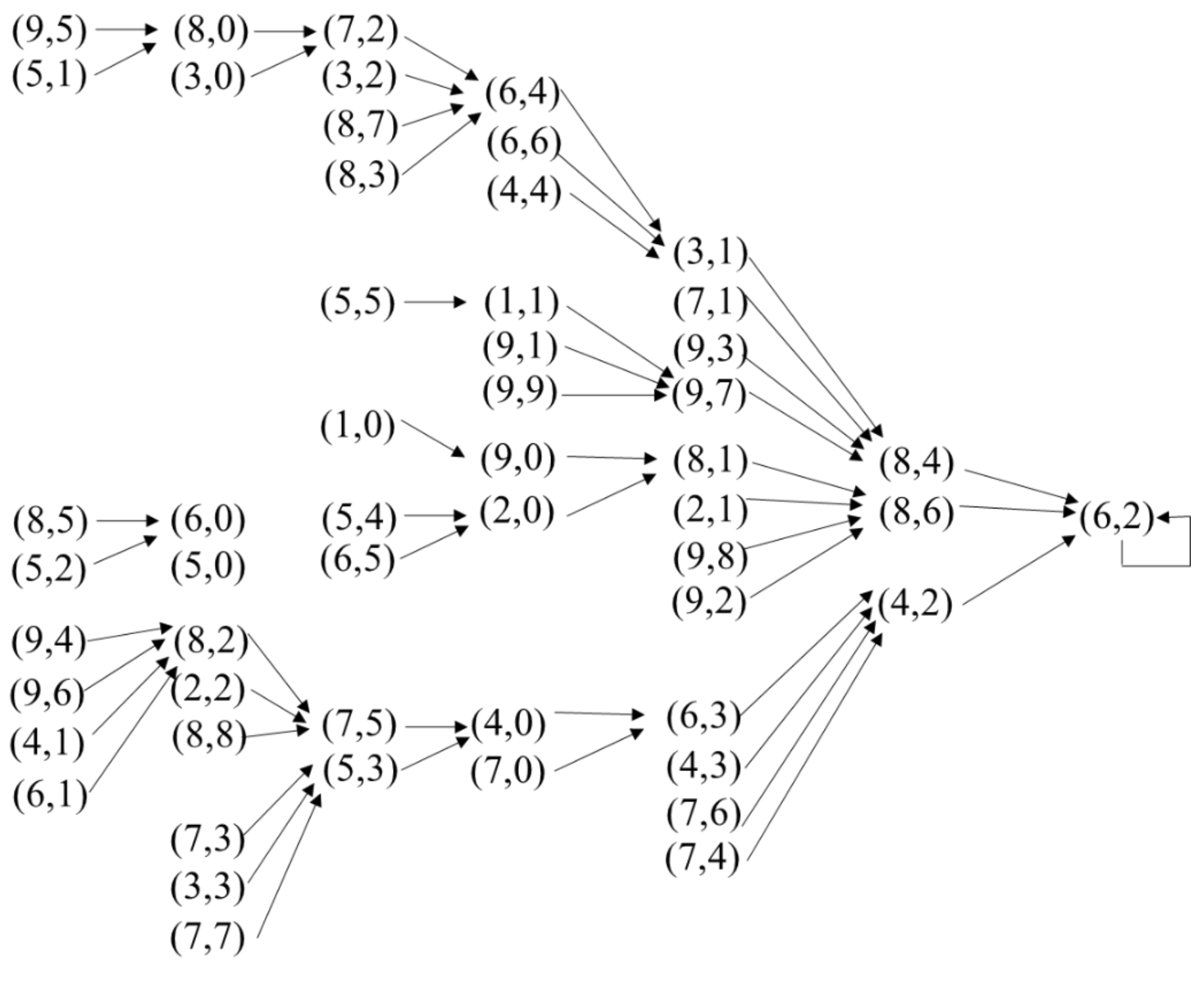

Graph 2. Tree of parametric transformations in $A_4$

The structure of the transformation tree can be analyzed by following several methodological approaches (Graph 3).

- Approach A: to study the relations between parameters shaping a particular layer. Each layer is defined by the distance or the number of transformations that separates it from the parameters of balance or the corresponding cycle. In $A_4(6,2)$.
- Approach B: to start from the closer layer to (6,2) and expand, encompassing a farther layer each time.
- Approach C: to study how the tree branches are formed (e.g., main branches, secondary branches), regardless of the layer or layers they may extend to.

  Now, we will develop approach C based on the concept of binary equivalence relation. This will allow us to establish equivalence classes corresponding to the different tree structures. The developed approach enables the definition of hierarchical partitions which finally converge with approach B.

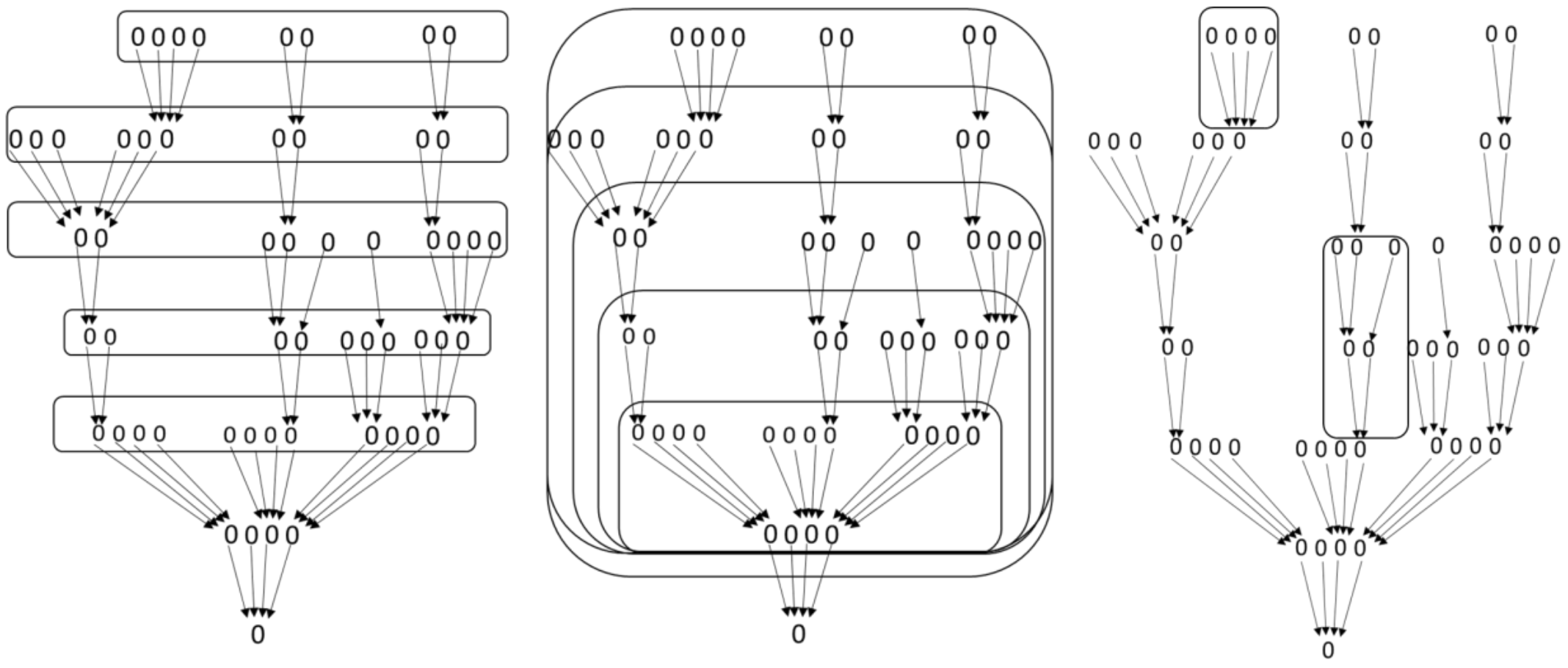

Approach A Approach B Approach C

Graph 3. Some approaches to the analysis of transformation trees

## 4. Equivalence relations

Two numbers m and n shall be considered equivalent $R_r$ of order r if and only if

$m\ R_r\ n \leftrightarrow K^r(m) = K^r(n) \qquad r \geq 0$ [1]

Since $K^r(m) = K^{r-1}[K(m)]$ the recurrence relations

$m\ R_r\ n \leftrightarrow K(m)\ R_{r-1}\ K(n) \qquad r \geq 1$ [2]

$m\ R_r\ n \leftrightarrow K^s(m)\ R_{r-s}\ K^s(n) \qquad s \leq r$ [3]

result.

Note that when two numbers are equivalent, they are so for any equivalence of a higher order

$m\ R_r\ n \rightarrow m\ R_{r+s}\ n \qquad s \geq 0$ [4]

It is necessary to distinguish **new** equivalences $\mathbf{R_r}$, which are those that are not of a lower order,

$m\ \mathbf{R_r}\ n \rightarrow m\ \not R_{r-1}\ n \qquad r > 0$ [5]

$m\ \mathbf{R_r}\ n \rightarrow m\ \not{R}_{r-s}\ n \quad s = 1,2\ldots r$ [5a]

from **old** equivalences ${}_aR_r$, which are the ones that are also of a lower order

$m\ ({}_aR_r)\ n \rightarrow m\ R_{r-1}\ n \qquad r > 0$ [6]

The notation $R_r$ does not specify whether the equivalence is new or old.

The equivalences of order r are binary equivalence relations that match the features

reflexive $m\ R_r\ m \qquad \forall\ m \in A_w$

symmetric $m\ R_r\ n \leftrightarrow n\ R_r\ m \qquad \forall\ m,n$

transitive $m\ R_r\ n$ and $n\ R_r\ q \rightarrow m\ R_r\ q \qquad \forall\ m,n,q$

They allow to group the numbers in equivalence classes of order r

$C_r = \{m,n \in A_w;\ m\ R_r\ n\}$ [7]

These equivalence classes are pairwise disjoint and their union is set $A_w$, generating a partition in $A_w$.

Due to the fact that when two numbers are equivalent they are so for any equivalence of a higher order [4], classes $C_r$ must include both new and old equivalences. This hierarchical and inclusive feature of binary relations $R_r$ is essential to understand the convergence of the process.

$\mathbf{C_r} = \{m,n \in A_w,\ m\ \mathbf{R_r}\ n\},\ {}_aC_r = \{\ p,q \in A_w,\ p\ R_{r-1}\ q\} \qquad C_r = \mathbf{C_r} \cup {}_aC_r$ [8]

The transitive feature and relation [4] yield

$m\ R_r\ n$ and $n\ R_{r+s}\ q \rightarrow m\ R_{r+s}\ q \qquad s \geq 0$ [9]

In our methodological approach to the analysis of the Kaprekar process it is essential to extend the concept of equivalence $R_r$ defined between numbers to sets of numbers with images common (equivalence classes).

Considering the last relationship, if

$m \in C_r^i,\ p \in C_s^j$ y $m \in R_t\ p,\ t = \max(r,s) + u,\ u \geq 0$ then if

$n \in C_r^i$ y $q \in C_s^j$ results $n\ R_t\ q$. This is due to the univocal character of the transformation. All numbers belonging to a class $C_v$ have the same image in the v-th transformation. From here on, any subsequent transformation will have the same image. Consequently, the element that best represents the transformation is not the number, but the equivalence class. Thus, we can agree to define the equivalence relation between classes according to

$m \in C_r^i,\ n \in C_s^j$ y $m\ R_t\ p,\ t = \max(r,s) + u,\ u \geq 0 \leftrightarrow C_r^i\ R_t\ C_s^j$ [10]

The direction of the arrow to the left must be understood to mean that if two classes $C_r^i$ y $C_s^j$ have a relation $R_t$ then all numbers m and n belonging to those classes have a relation $R_t$. In particular turns out

$m \in C_r^i,\ \ n \in C_r^j,\ \ m\ R_{r+s}\ n \leftrightarrow C_r^i\ R_{r+s}\ C_r^j \quad s \geq 0$ [10a]

Due to basic functions f, numbers with the same parameters α and ß yield the same image and only them. As we indicated in section 2,

$p(m) = p(n) \leftrightarrow K(m) = K(n)$

Consequently, taking into account [1] it results

$m\ R_r\ n \leftrightarrow p\ [K^{r-1}(m)] = p\ [K^{r-1}(n)] \qquad r \geq 1$ [11]

## 5. Equivalences $R_0$ and $R_1$

a) $R_0$ is the identity relation. From [1]

$m\ R_0\ n \leftrightarrow K^0(m) = K^0(n) \leftrightarrow m = n$

each equivalence class $C_0$ is formed by a single number. There are thus 9990 classes $C_0$ in $A_4$

b) $R_1$ links numbers which yield the same image with only one transformation. From [1]

$m\ R_1\ n \leftrightarrow K(m) = K(n)$ and taking [11] into account

$m\ R_1\ n \leftrightarrow p(m) = p(n)$

Equivalence classes $C_1$ will be those consisting of all numbers which have the same parameters α and ß. For simplicity purposes, classes $C_1$ shall be represented as (α,ß), since this way they are characterized:

$(\alpha,ß) = \{ n \in A_w,\ p(n) = (\alpha,ß)\}$

The double use of (α,ß), both as parameters of a number and as a set of numbers with the same parameters is not ambiguous in the context in use. In $A_4$ and $A_5$ cases, classes $C_1$ are defined by two parameters; by one parameter in $A_2$ and $A_3$ and by 3 in $A_6$ and $A_7$. We use binary representation in a general way for convenience.

Classes $C_1 = (\alpha,ß)$ play a key role in this paper, as we have developed functions capable of transforming them: functions $f_i$ provide the numerical image of a class, while functions $K_i$ (Tables 2 and 3 in the Annex), $K_i(\alpha,ß) = (\alpha', ß')$, provide the parameters of the transformed number. They have a double use:

- On the one hand, as noted at the beginning of section 3, they simplify the transformation process by switching from numbers to classes (from 9990 to 54 in $A_4$)
- On the other hand, along with functions $K_i$ which operate on them, they allow to build higher-order equivalences

This is simplified because basic equivalence equations can be translated thanks to [10] and [11] in terms of class. Thus,

$m \in (\alpha_1, ß_1),\ n \in (\alpha_2, ß_2),\ m' = K(m) \in (\alpha'_1, ß'_1),\ n' = K(n) \in (\alpha'_2, ß'_2)$

taking [2] into account

$(\alpha_1, ß_1)\ R_r\ (\alpha_2, ß_2) \leftrightarrow (\alpha'_1, ß'_1)\ R_{r-1}\ (\alpha'_2, ß'_2) \qquad r \geq 1$ [12]

## 6. Approaches to find equivalences

[7] and [9] yield

$(\alpha_1, ß_1) \subset C^i_r,\ (\alpha_2, ß_2) \subset C^j_r,\ (\alpha_1, ß_1)\ R_{r+s}\ (\alpha_2, ß_2) \rightarrow C^i_r\ R_{r+s}\ C^i_r$ [13]

To find equivalences of order $r>1$ it suffices to find them between classes (α,ß).
A function or operator $e_r$ can be linked to any equivalence $R_r$, such that

$(\alpha_1,ß_1)\ R_r\ (\alpha_2,ß_2) \leftrightarrow e_r\ (\alpha_1,ß_1) = (\alpha_2,ß_2)$ [14]

It is a type of representation that makes calculation easier. The product of functions or operators will be used as usually

$(f \times g)\ (\alpha,ß) = f\ [g(\alpha,ß)]$

For example in $A_4$,

$K_{17}\ (\alpha,ß) = [(\alpha-ß)+1, (\alpha-ß)-1]$, $e_2^5\ (\alpha,ß) = (10-ß, 10-\alpha)$

$(e_2^5 \times K_{17})\ (\alpha,ß) = e_2^5\ [K_{17}\ (\alpha,ß)] = [11-(\alpha-ß), 9-(\alpha-ß)]$

whose validity is limited to $K_{17}$ and $e_2^5$ 's domains of existence.

Approach A: Identifying sequences

If two classes (α,ß) generate numerical images which are permutations of the same digits, then such classes will be equivalent $R_2$

$f\ (\alpha_1,ß_1) = m$, $f\ (\alpha_2,ß_2) = n$, $m=P(n) \rightarrow (\alpha_1,ß_1)\ R_2\ (\alpha_2,ß_2)$

as m and n will yield the same sequence when sorting their digits to apply Kaprekar's process, and therefore the image of such sequence will be the same.

This approach is used to find equivalence $R_2^1$ and equivalences $R_2$ in group I, $R_2^2$, $R_2^3$, $R_2^4$ in $A_4$.

Approach B: Identifying classes

It is based on relation [12]. Specifically, with r=2

$(\alpha_1,ß_1)\ R_2\ (\alpha_2,ß_2) \leftrightarrow (\alpha'_1,ß'_1) = (\alpha'_2,ß'_2)$

Let us transform classes $(\alpha_1,ß_1)$ and $(\alpha_2,ß_2)$ with functions $K_i$ and impose the identity between their images

This is the approach used to establish the equivalences $R_2$ of set II: $R_2^5$, $R_2^6$, $R_2^7$, $R_2^8$ in $A_4$.

Approach C: Hierarchization of relations

Uses relation [12] directly

$(\alpha_1,ß_1)\ R_r\ (\alpha_2,ß_2) \leftrightarrow (\alpha'_1,ß'_1)\ R_{r-1}(\alpha'_2,ß'_2) \quad r\geq 1$

Operationally, the basic condition is

$e_{r-1}\ [K_i\ (\alpha_1,ß_1)] = K_j\ (\alpha_2,ß_2)$

It demands using not only functions $K_i$, as in approach B, but also known $R_{r-1}$.

This approach is therefore useful to find equivalences of an increasingly higher order.

It has been used to establish equivalences $R_3$ and of a higher order.

Approach D: Obtaining inverse images

Uses relation [12] in sense to left

$(\alpha_1,ß_1)\ R_{r-1}\ (\alpha_2,ß_2) \rightarrow K^{-1}\ (\alpha_1,ß_1)\ R_r\ K^{-1}(\alpha_2,ß_2)$

Its validity is subject to the classes connected by $R_{r-1}$ having inverse images. If there is a solution, it does not have to be unique, as there may be several anti-images.

For example, classes (9,5) and (5,1) do not have inverse images. (9,5) $R_2^5$ (5,1) but the corresponding transformation chains get cut and there is no $R_3$. In most cases of classes $C_2$ only one class (α,ß) has inverse images, while the rest do not.

In $A_4$ only 19 out of 54 existing classes (α,ß) have inverse images.

As we increase the order in the hierarchy, the new equivalence relations only occur between a very reduced number of classes. Often only between one pair of classes. In these situations, approach D is very useful.

Approach E: Product of equivalences

From known equivalences, secondary equivalences can be obtained through the product of compatible equivalences. For example, in $A_4$

$e_2^2$ (α,ß) = (10-α, ß) $\alpha + ß \leq 10$

$e_3^{21}$(α, 9-α) = (14-α, α-5) $5 < \alpha < 9$

$e_3^{211} = e_2^2 \times e_3^{21}$ (α,9- α) = (α-4, α-5) $5 < \alpha < 9$

Thus, (8,1) $e_3^{211}$ (4,3)

The order of the new equivalence is determined by the order of the higher-order operator, as requires [9]

Approach F: Using operators $K^r$, $r > 1$

It is based on [3], which when translated to classes can be written

$K^r$ (α1,ß1) $R_s$ $K^r$ (α2,ß2) → (α1,ß1) $R_{s+r}$ (α2,ß2) $r > 1$

The basic condition is

$e_s[K^r$ (α1, ß1)] = $K^r$ (α2,ß2)

In order to obtain $K^r$, $K_i$ that are compatible must be used (Tables 1 and 2 in the Annex). The second operator must have in its domain of existence the image of the first, and so on. For instance, in $A_4$ $K_{26}$ can act on 10 functions $K_i$, while $K_{25}$ acts on 5 and $K_9$ only on $K_{17}$. Add to this the compatibility with the equivalence function $e_s$.

Let the diagram of a process r=2, s=2 in $A_4$ be as follows

$$\left.\begin{array}{l} (9,0) \xrightarrow{K_{25}} (8,1) \xrightarrow{K_{17}} (8,6) \\ (4,0) \xrightarrow{K_{26}} (6,3) \xrightarrow{K_{17}} (4,2) \end{array}\right] \quad (8,6)\ R_2^5(4,2) \rightarrow (9,0)\ R_4\ (4,0)$$

$$\left.\begin{array}{l} K^2\ (9,0) = K_{17} \times K_{25}\ (9,0) = (8,6) \\ K^2\ (4,0) = K_{17} \times K_{26}\ (4,0) = (4,2) \end{array}\right] \quad (9,0)\ R_{2+2}\ (4,0)$$

The approach becomes the most useful when used in sequence (section 7.4). As r increases, the domains of existence of functions $K^r$ decrease, sometimes being limited to a single class.

In general, the approaches mentioned above demand that

- The used functions $K_i$ be compatible with the classes (α,ß) to be transformed
- The equivalence $e_r$ be applied on a class within its domain of existence
- All classes have an image, but the reciprocal is not always true. When working with inverse images, its existence must be ensured.
-

## 7. Equivalences in $A_4$

Let us recall that

$f_1(\alpha,ß) = (\alpha \;\; ß\text{-}1 \;\; 9\text{-}ß \;\; 10\text{-}\alpha) \;\; \alpha \geq ß, \;\; 0 < ß \leq 9$

$f_2(\alpha,0) = (\alpha\text{-}1 \;\; 9 \;\; 9 \;\; 10\text{-}\alpha) \;\; 0 < \alpha \leq 9$

### 7.1 Equivalences $R_1$

Taking into account α and ß's definition, the corresponding classes include all permutations $P_i$ of those numbers that satisfy

$(\alpha,ß) = \{P_i\,(\alpha+d \;\; ß+c \;\; c \;\; d;\; c \geq d,\; \alpha+d \geq ß+c\} \;\forall_i = 1, 2,\ldots 24$

In Graph 1, 1998 $R_1$ 8991, numbers which are permutations of the same digits and both belonging to class (8,1). Also 1904 $R_1$ 2509 with different digits but both belonging to class (9,3).

It is worth noting class (α,0), if

$n \in (\alpha,0) \rightarrow O(n) = (a\, b\, b\, d),\;\; a \geq b, b \geq d,\;\; \alpha = a\text{-}d > 0$

For each set a, b and d there are 10 different equivalent permutations. For example, let (a b b d) $R_1$ (d a b b)

Since $0 < \alpha \leq 9$, $0 \leq ß \leq 9$ and $\alpha \geq ß$ there are $\sum_{s=2}^{10} s = 54$ different classes $C_1$.

Sometimes it is useful to remember that number n =(α ß 0 0) belongs to class (α,ß) and therefore if

$m \in (\alpha,ß) \rightarrow K(m) = K(\alpha \;\; ß \;\; 0 \;\; 0) \leftrightarrow m\, R_1\, (\alpha \;\; ß \;\; 0 \;\; 0)$

The transformation tree between classes $C_1$ = (α,ß) appears in Graph 2.

### 7.2 Equivalences $R_2$

All equivalences $R_2$ can be reduced to equivalences $R_1$ between numbers belonging to $B_4$, since due to [2] and $B_4$'s definition

$m\, R_2\, n \leftrightarrow K(m) \;\; R_1 \;\; K(n)$ , $K(m) \in B_4$, $K(n) \in B_4$

**Group ß = 0** (Approach A)

Due to $f_2$, in $B_4$ there can only be two permutations with the same digits of a number $P_{12}$ = (a 9 9 9-a) and $P_{21}$ = (9-a 9 9 a) $a \leq 8$ so $K(P_{12}) = K(P_{21})$ and the classes that generate $P_{12}$ and $P_{21}$ will be $R_2$.

That is,

m ∈ ($\alpha_1$,0), n ∈ ($\alpha_2$,0), K(m) = ($\alpha_1$-1 9 9 10-$\alpha_1$), K (n) = ($\alpha_2$-1 9 9 10-$\alpha_2$)

- **Equivalence $R_2^0$ ($\alpha_1$,0) = ($\alpha_2$,0)**
  K(n) = $P_{12}$ [K(m)] ↔ ($\alpha_1$-1 9 9 10-$\alpha_1$) = ($\alpha_2$-1 9 9 10-$\alpha_2$) ↔ $\alpha_1 = \alpha_2$
  It is the identity equivalence: $R_2^0 = R_1$
- **Equivalence $R_2^1$ $\alpha_1+\alpha_2=11$**
  K(n) = $P_{21}$ [K(m)] ↔ (10-$\alpha_1$ 9 9 $\alpha_1$-1) = ($\alpha_2$-1 9 9 10-$\alpha_2$) ↔10-$\alpha_1$ = $\alpha_2$-1
  Thus, in Graph 1, 0999 $R_2$ 4446
  m = 0999 ∈ (9,0), n= 4446 ∈ (2,0) $\alpha_1+\alpha_2$=9+2 = 11
  m' = K(m) = 8991∈ (8,1), n' = K(n) = 1998 = $P_{21}$ [K(m)] ∈ (8,1)
  m" = $K^2$ (m) = 8082 = $K^2$(n)
  This equivalence groups classes $C_1$ = (α,ß) into five classes $C_2$

Table 1. Classes $C_2$ of equivalence $R_2$ with ß = 0 in $A_4$

| Class $C_2$ | $C_2^{11}$ | $C_2^{7}$ | $C_2^{14}$ | $C_2^{8}$ | $C_2^{15}$ |
|---|---|---|---|---|---|
| Classes $C_1$ forming it | (1,0) | (2,0)<br>(9,0) | (3,0)<br>(8,0) | (4,0)<br>(7,0) | (5,0)<br>(6,0) |
| Image (α', ß') ∈ $B_4$ | (9,0) | (8,1) | (7,2) | (6,3) | (5,4) |

An equivalence function can be associated with each equivalence, as follows
$R_2^0$ : $e_2^0$ (α,0) = (α, 0) ; $R_2^1$ : $e_2^1$ (α,0) = (11-α,0), α > 1
$e_2^1$ is an involution, i.e., it matches its own inverse $e_2^1 = [e_2^1]^{-1}$

The product of these functions $e_2^i$ x $e_2^j$ has an isomorphic abelian group algebraic structure to $\mathbb{Z}_2$ and the symmetric group of two permutations $S_2$.

**<u>Groups ß>0</u>**

In $B_4$ there can only be four permutations having the same digits of a number
$P_{1234}$ = (a b 8-b 10-a), $P_{1324}$ = (a 8-b b 10-a)
$P_{4231}$ = (10-a b 8-b a), $P_{4321}$ = (10-a 8-b b a), a≥b+1

since only those satisfy the conditions imposed by $f_i$ for images in classes (α,ß) with ß>0. For the same reasons mentioned for group ß = 0, the classes generated will be equivalent $R_2$.

It is important not to confuse the possible permutations of n' = K (n) aforementioned with permutations O (n') mentioned in Tables 2 and 4 in the Annex. For example, O (n') cannot have permutations $P_{4231}$ and $P_{4321}$, but n' can.

There are other permutations of n’, but they do not belong to $B_4$, although they do belong to the same class $C_1$ as the others do, and will yield the same image. Accordingly, there is no class $C_1$ completely included in $B_4$, and all of them have their image in it.

The possible permutations rearrange pairs of elements:

Let us take $P_{1234}$ as a base sequence. No transposition of elements.

$P_{4231}$ is a transposition of the extreme digits but keeps the middle ones in place. Therefore it replaces a with 10-a and vice versa.

$P_{1324}$ is a transposition of the middle digits but keeps the extreme ones in place. It replaces b with 8-b and vice versa.

$P_{4312}$ rearranges both pairs. It is the product of the two previous transpositions.

There exist two basic sequences which are compatible with $B_4$

$S_1 = (\alpha \;\; ß\text{-}1 \;\; 9\text{-}ß \;\; 10\text{-}\alpha)$ resulting from $a = \alpha$, $b = ß\text{-}1$

$S_2 = (ß \;\; \alpha\text{-}1 \;\; 9\text{-}\alpha \;\; 10\text{-}ß)$, $a = ß$ and $b = \alpha\text{-}1$

Both replace α with ß and vice versa. Each one of them can arise from different configurations of parameters α and ß.

For example, if $p(n) = (10\text{-}ß, 10\text{-}\alpha)$, then $n' = K(n) = (10\text{-}ß \;\; 9\text{-}\alpha \;\; \alpha\text{-}1 \;\; ß) = P_{4321}$ (S2) and if $p(m) = (10\text{-}\alpha, ß)$ then $m' = (10\text{-}\alpha \;\; ß\text{-}1 \;\; 9\text{-}ß \;\; \alpha) = P_{4231}$ ($S_1$)

The possible permutations need to be images of different classes $C_1$, with the natural exception of the identity, as each class yields a unique image. But what are those classes?

Equivalences within the same sequence: Group I

For the sake of brevity, let us use the base sequence $S_1$. The same types appear between rearrangements of sequence $S_2$.

Let us use approach A. Remember that

$m \in (\alpha_1, ß_1)$, $m' = K(m) = (\alpha_1 \;\; ß_1\text{-}1 \;\; 9\text{-}ß_1 \;\; 10\text{-}\alpha_1)$ $\;\; ß_1 > 0$

If a certain permutation P of the image of

$n \in (\alpha_2, ß_2)$, $n' = K(n) = (\alpha_2 \;\; ß_2\text{-}1 \;\; 9\text{-}ß_2 \;\; 10\text{-}\alpha_2)$ $\;\; ß_2 > 0$

matches m’

$P(n') = m'$, then $p(n') = p(m') \leftrightarrow K^2(m) = K^2(n) \leftrightarrow m \, R_2 \, n$

**Equivalence $\mathbf{R_2^0}$**: $(\alpha_1, ß_1) = (\alpha_2, ß_2) \leftrightarrow m \, R_1 \, n \leftrightarrow m \, R_2 \, n$, Identity

**Equivalence $\mathbf{R_2^2}$**: $\alpha_2 = 10\text{-}\alpha_1$, $ß_2 = ß_1$ $\;\; \alpha_1 + ß_1 \leq 10$, $\alpha_1 \geq ß_2$

It results from $P_{1234}(m') = P_{4231}(n') = (10\text{-}\alpha_2 \;\; ß_2\text{-}1 \;\; 9\text{-}ß_2 \;\; \alpha_2)$

$\alpha_1 = 10-\alpha_2$, $ß_1-1 = ß_2-1$, $9-ß_1 = 9-ß_2$, $10-\alpha_1 = \alpha_2 \leftrightarrow \alpha_1 + \alpha_2 = 10$, $ß_1 = ß_2 \leftrightarrow m\ R_2^2\ n$

Since it is a transposition of the extreme digits, this can be simplified by replacing α with 10-α

$(\alpha,ß)\ R_2^2\ (10-\alpha,ß)$

**Equivalence $R_2^3$:** $\alpha_2 = \alpha_1$ $ß_2 = 10-ß_1$ $\alpha_1+ß_1 \geq 10$, $\alpha_1 \geq ß_1$

It results from $P_{1234}$ (m') = $P_{1324}$ (n') : $\alpha_1 = \alpha_2$, $ß_1-1 = 9-ß_2$
$(\alpha,ß)\ R_2^3\ (\alpha,10-ß)$

Thus, in Graph 1, (m=9261) $R_2^2$ (n=8082) since p (m) = (8,4), p (n) = (8,6), $\alpha_1 = \alpha_2 = 8$, $ß_1+ß_2 = 4+6 = 10$ and it is verified that $\alpha_1 + ß_1 = 8+4 \geq 10$ then K (m) = 8352 and K (n) = 8532 = $P_{1324}$(8352) p [K(m)] = p[K(n)] = (6,2) and $K^2$ (m) = $K^2$ (n) = 6174

**Equivalence $R_2^4$** $\alpha_2 = 10-\alpha_1$, $ß_2 = 10-ß_1$, $\alpha_1 = ß_1$

It results from $P_{1234}$ (m) = $P_{4321}$ (n') ; $\alpha_1 = 10-\alpha_2$, $ß_1-1 = 9-ß_2$

Associated equivalence functions

$R_2^0$: $e_2^0(\alpha,ß) = (\alpha,ß)$; $R_2^2$: $e_2^2(\alpha,ß) = (10-\alpha,ß)$

$R_2^3$: $e_2^3(\alpha,ß) = (\alpha,10-ß)$; $R_2^4$: $e_2^4(\alpha,ß) = (10-\alpha,10-ß)$

The product of these functions is commutative and it provides its set with an isomorphic abelian group structure of Klein group or $D_2$.

Table 2. Klein group of equivalences $R_2$ in $A_4$

(g x f) (α,ß) = g [f(α,ß)]

| g \ f | $e_2^0$ | $e_2^2$ | $e_2^3$ | $e_2^4$ |
|---|---|---|---|---|
| $e_2^0$ | $e_2^0$ | $e_2^2$ | $e_2^3$ | $e_2^4$ |
| $e_2^2$ | $e_2^2$ | $e_2^0$ | $e_2^4$ | $e_2^3$ |
| $e_2^3$ | $e_2^3$ | $e_2^4$ | $e_2^0$ | $e_2^2$ |
| $e_2^4$ | $e_2^4$ | $e_2^3$ | $e_2^2$ | $e_2^0$ |

Equivalences appear as Klein subgroups in the symmetric group $S_4$ of all possible permutations of 4 elements.

Equivalence between sequences: Set II

All equivalences m $R_2$ n in group I generate images m' and n' which are permutations of the same sequence of digits. In the previous paragraph we referred to base sequence $S_1$ (α ß-1 9-ß 10-α) to simplify. However, all results are valid for sequence $S_2$ (ß α-1 9-α 10-ß) just by replacing α with ß and vice versa. Only these two sequences – with different digits– and the transposition of their extreme and/or middle digits are compatible with $B_4$.

The replacement α ↔ ß automatically generates four mutually equivalent classes.

(α,ß) $R_2$ (10-α, ß) $R_2$ (α,10-ß) $R_2$ (10-α, 10-ß) ↔

(ß,α) $R_2$ (10-ß, α) $R_2$ (ß,10-α) $R_2$ (10-ß, 10-α)

The images of these classes correspond to the 4 possible permutations of $S_2$

$(ß,\alpha) \xrightarrow{K} S_2 = P_{1234}(S_2)$; $(10-ß, \alpha) \rightarrow P_{4231}(S_2)$;

$(ß, 10-\alpha) \rightarrow P_{1324}(S2)$ ; $(10-ß, 10-\alpha) \rightarrow P_{4321}(S_2)$

And the associated functions will be identical, already mentioned in group I

$e_2^0$ (α,ß) = (ß,α); $e_2^2$ (ß,α) = (10-ß, α); $e_2^3$ (ß,α)= (ß, 10-α); $e_2^4$ (ß,α) = (10-ß,10-α)

It would just be a matter of letter shift.

But is there any equivalence relation between numbers of the two sequences?

m = 1998 ∈ (8,1), m' = 8082 ∈ (8,6), m' = $P_{1234}(S_1)$

n = 3095 ∈ (9,2), n' = 9171 ∈ (8,6), n' = $P_{4321}(S_2)$

$\alpha_1 = 8$, $ß_1 = 1$, $\alpha_2 = 10-ß_1 = 9$, $ß_2 = 10-\alpha_1 = 2$

Now it is not just a matter of letter shift. Is there an equivalence relation m $R_2$ n? The necessary and sufficient condition is established by the relation [12], approach B

(9,2) $R_2$ (8,1) ↔ ($\alpha'_1$, $ß'_1$) = ($\alpha'_2$, $ß'_2$)

and in this case $\alpha'_1 = \alpha'_2 = 8$, $ß'_1 = ß'_2 = 6$ is verified

Is there any function $K_i$ which yields images (α', ß') that do not vary when replacing α with ß and vice versa? Initially, this is what happens with the four functions $K_i$ belonging to type 2 (Table 2 in the Annex). Functions $K_9$, $K_{10}$, $K_{13}$ and $K_{14}$ yield images that depend on the addition. For example,

$K_9$ (α,ß) = [(α+ß)-9, (α+ß)-11], α ≤ ß+1, α+ß ≥ 11

Replacing α with ß and vice versa would not alter the addition. But the strict conditions to be satisfied by α and ß so $K_9$ can operate mean that there are no two classes with an equal sum that satisfy such conditions. The same is true for $K_{10}$, $K_{13}$ and $K_{14}$.

**Equivalence $R_2^5$:** $\alpha_2 = 10-ß_1$, $ß_2 = 10-\alpha_1$ $\alpha_1 \geq ß_1$

However, there are other classes $(\alpha_1,ß_1)$ and $(\alpha_2,ß_2)$ which generate permutations of $S_1$ and $S_2$ and are compatible with functions $K_i$ of type 3 (Table 2 in the Annex). Such is the case of $(\alpha_1,ß_1)$ and $(10-ß_1, 10-\alpha_1)$ with functions $K_{17}$ and $K_{18}$. Thus,

$K_{17}(\alpha,ß) = [(\alpha-ß)+1, (\alpha-ß)-1]$ $\alpha \geq ß+1$, $9 \leq \alpha+ß \leq 11$

If $\alpha_1-ß_1 = \alpha_2-ß_2 = (10-ß_1) - (10-\alpha_1) \rightarrow K_{17}(\alpha_1,ß_1) = K_{17}(\alpha_2,ß_2) \leftrightarrow (\alpha_1,ß_1)\ R_2\ (\alpha_2,ß_2)$

It is the case of classes (8,1) and (9,2), referred to above. Thus arises the equivalence $(\alpha_1,ß_1)\ R_2^5\ (10-ß_1,10-\alpha_1)$ which relates sequences $S_1$ and $S_2$.

**Equivalence $R_2^6$:** $\alpha_2 = ß_1$, $ß_2 = 10-\alpha_1$ $\alpha_1+ß_1 \geq 10$

Having established equivalence $(\alpha,ß)\ R_2^5\ (10-ß, 10-\alpha)$, due to the transitive feature of equivalences, equivalences between (α,ß) –in addition to the other equivalences in group I– and those associated with the transpositions in $S_2$ can be automatically established.

Thus, $(\alpha_1,ß_1)\ R_2^6\ (ß_1, 10-\alpha_1)$

It can be also interpreted as the result of applying equivalence $R_2^2$ to $R_2^5$

**Equivalence $R_2^7$:** $\alpha_2 = 10-ß_1$, $ß_2 = \alpha_1$, $\alpha_1+ß_1 \leq 10$

$(\alpha_1,ß_1)\ R_2^7\ (10-ß_1, \alpha_1)$; $R_2^3 \times R_2^5 = R_2^7$

**Equivalence $R_2^8$** $\alpha_2 = ß_1$, $ß_2 = \alpha_1$, $\alpha_1 = ß_1$

$(\alpha_1,ß_1)\ R_2^8\ (ß_1,\alpha_1)$; $R_2^4 \times R_2^5 = R_2^8$

If $ß_1 = 10-\alpha_1$, $R_2^5 = R_2^0$, $R_2^6 = R_2^2$, $R_2^7 = R_2^3$, $R_2^8 = R_2^4$

<u>Associated equivalence functions</u>

$R_2^5$: $e_2^5(\alpha,ß) = (10-ß, 10-\alpha)$; $R_2^6$: $e_2^6(\alpha,ß) = (ß, 10-\alpha)$

$R_2^7$: $e_2^7(\alpha,ß) = (10-ß, \alpha)$; $R_2^8$: $e_2^8(\alpha,ß) = (ß,\alpha)$

The domains of existence of all functions $R_2$ are shown in Table 3.

The product of these functions of set II is not a closed binary operation on the set. All products yield functions from group I, including $e_2^0$ which can be assigned to both sets. Equivalences $e_2^5$ and $e_2^8$ are involutory; $e_2^6$ and $e_2^7$ are not.

The product of functions R2 is shown in the table below (Table 4). It has a non-abelian group structure where the Klein four-group $\{e_2^0, e_2^2, e_2^3, e_2^4\}$ is not a normal subgroup.

Table 3. Equivalence functions R2 in A4

| | | |
|---|---|---|
| ß=0 | α = 1 | $e_2^0$ (α,0) = (α,0) = (1,0) |
| | α> 1 | $e_2^1$ (α,0) = (11-α,0) |
| ß>0 | α = ß | $e_2^4$ (α,ß) = (10-α,10-ß) |
| | | $e_2^8$ (α,ß) = (ß,α) |
| | α ≥ ß | $e_2^0$ (α,ß) = (α,ß) |
| | | $e_2^5$ (α,ß) = (10-ß ,10-α) |
| | $\alpha + ß \geq 10$ | $e_2^3$ (α,ß) = (α,10-ß) |
| | | $e_2^6$ (α,ß) = (ß ,10-α) |
| | $\alpha + ß \leq 10$ | $e_2^2$ (α,ß) = (10-α,ß) |
| | | $e_2^7$ (α,ß) = (10-ß,α) |

Table 4. Group of equivalences R2 in A4: (gxf)(α,ß) = g [f(α,ß)]

| g \ f | $e_2^0$ | $e_2^2$ | $e_2^3$ | $e_2^4$ | $e_2^5$ | $e_2^6$ | $e_2^7$ | $e_2^8$ |
|---|---|---|---|---|---|---|---|---|
| $e_2^0$ | $e_2^0$ | $e_2^2$ | $e_2^3$ | $e_2^4$ | $e_2^5$ | $e_2^6$ | $e_2^7$ | $e_2^8$ |
| $e_2^2$ | $e_2^2$ | $e_2^0$ | $e_2^4$ | $e_2^3$ | $e_2^6$ | $e_2^5$ | $e_2^8$ | $e_2^7$ |
| $e_2^3$ | $e_2^3$ | $e_2^4$ | $e_2^0$ | $e_2^2$ | $e_2^7$ | $e_2^8$ | $e_2^5$ | $e_2^6$ |
| $e_2^4$ | $e_2^4$ | $e_2^3$ | $e_2^2$ | $e_2^0$ | $e_2^8$ | $e_2^7$ | $e_2^6$ | $e_2^5$ |
| $e_2^5$ | $e_2^5$ | $e_2^7$ | $e_2^6$ | $e_2^8$ | $e_2^0$ | $e_2^3$ | $e_2^2$ | $e_2^4$ |
| $e_2^6$ | $e_2^6$ | $e_2^8$ | $e_2^5$ | $e_2^7$ | $e_2^2$ | $e_2^4$ | $e_2^0$ | $e_2^3$ |
| $e_2^7$ | $e_2^7$ | $e_2^5$ | $e_2^8$ | $e_2^6$ | $e_2^3$ | $e_2^0$ | $e_2^4$ | $e_2^2$ |
| $e_2^8$ | $e_2^8$ | $e_2^6$ | $e_2^7$ | $e_2^5$ | $e_2^4$ | $e_2^2$ | $e_2^3$ | $e_2^0$ |

In set II symmetry is broken. $R_2^5$ and $R_2^8$ are symmetric, but $R_2^6$ and $R_2^7$ are not.

$(\alpha_1,ß_1)\ R_2^6\ (\alpha_2,ß_2) \leftrightarrow (\alpha_2,ß_2)\ R_2^7\ (\alpha_1,ß_1)$

Or, equivalently

$e_2^6 \text{ x } e_2^6 = e_2^4 \neq e_2^0$, $e_2^7 \text{ x } e_2^7 = e_2^4 \neq e_2^0$, $e_2^6 \text{ x } e_2^7 = e_2^0$, $e_2^7 \text{ x } e_2^6 = e_2^0$

However, the equivalence of second order stays

$(\alpha_1,ß_1)\ R_2\ (\alpha_2,ß_2) \leftrightarrow (\alpha_2,ß_2)\ R_2\ (\alpha_1,ß_1)$

as does the other essential feature, transitivity. This is why equivalence set II is important for establishing class groupings.

In addition to the 5 equivalence classes $R_2$ with ß=0 (Table 1), equivalences $R_2^2$ to $R_2^8$ with ß>0 generate 5 more classes that are shown in Table 5.

Table 5. Equivalence classes $C_2$ with ß > 0 in $A_4$

| Equivalence classes $C_2$ | $C_2^1$ | $C_2^2$ | $C_2^3$ | $C_2^4$ | $C_2^5$ | $C_2^6$ | $C_2^9$ | $C_2^{10}$ |
|---|---|---|---|---|---|---|---|---|
| Classes $C_1$=(α,ß) included | (8,6)<br>(8,4)<br>(6,2)<br>(4,2) | (9,7)<br>(9,3)<br>(7,1)<br>(3,1) | (9,8)<br>(9,2)<br>(8,1)<br>(2,1) | (7,6)<br>(7,4)<br>(6,3)<br>(4,3) | (6,6)<br>(6,4)<br>(4,4) | (9,9)<br>(9,1)<br>(1,1) | (8,7)<br>(8,3)<br>(7,2)<br>(3,2) | (5,5) |
| Image of $C_2$ | (6,2) | (8,4) | (8,6) | (4,2) | (3,1) | (9,7) | (6,4) | (1,1) |

Table 5 (follow-up)

| $C_2^{12}$ | $C_2^{13}$ | $C_2^{16}$ | $C_2^{17}$ | $C_2^{18}$ | $C_2^{19}$ | $C_2^{20}$ |
|---|---|---|---|---|---|---|
| (6,5)<br>(5,4) | (7,5)<br>(5,3) | (8,8)<br>(8,2)<br>(2,2) | (7,7)<br>(7,3)<br>(3,3) | (9,5)<br>(5,1) | (8,5)<br>(5,2) | (9,6)<br>(9,4)<br>(6,1)<br>(4,1) |
| (2,0) | (4,0) | (7,5) | (7,5) | (8,0) | (6,0) | (8,2) |

With these tables 1 and 5, Graph 4 can be built (upper row), which shows the structure of the convergence process of classes $C_2$ towards class (6,2). The convergence level states the number of transformations necessary to become class (6,2). For instance, numbers in class (3,0) $\in C_2^{14}$ such as 7545, and those in class (8,2) $\in C_2^{16}$ such as 2914 require 5 transformations for their successive images to belong to class (6,2), and 6 transformations to become 6174.

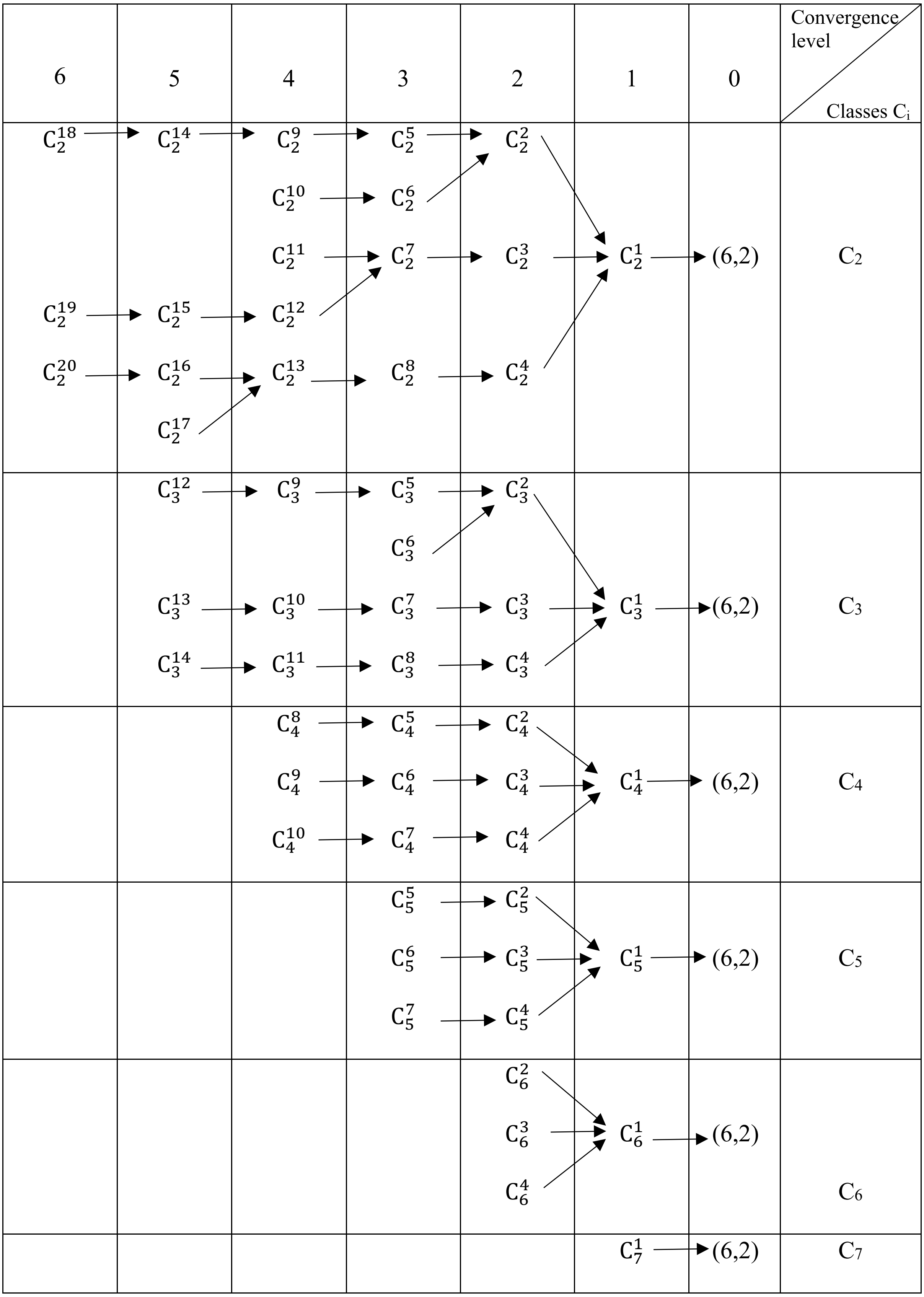

Graph 4. Transformation relations between the equivalence classes $C_2$ to $C_7$ in $A_4$

To delve further into the convergence process it is necessary to know the structure of equivalences $R_3$ and above.

### 7.3 Equivalences $R_3$

New equivalences **$R_3$** – those not being mere $R_2$, since all $R_2$ are $R_3$ – are established between classes whose images are related by new **$R_2$** [5]. Let us use approach C, [12].

$(\alpha_1,ß_1)\ R_3\ (\alpha_2,ß_2) \leftrightarrow (\alpha'_1,ß'_1)\ R_2\ (\alpha'_2,ß'_2)$, with the condition

$e_2 \times K_i\ (\alpha_1,ß_1) = K_j\ (\alpha_2,ß_2)$

**a) $K_i = K_j$**

This happens between the images of all functions $K_i$, with the exception of those operating on only one class ($K_{21}$), two classes ($K_{10}$ and $K_{13}$) or three classes ($K_{18}$) or $K_{25}$ and $K_{26}$ which operate on classes $(\alpha,0)$.

**Equivalences $R_3^1$**

For instance, $K_{14}$ generates classes with equivalence $R_2^5$: $e_2^5\ (\alpha,ß) = (10\text{-}ß, 10\text{-}\alpha)$

$K_{14}\ (\alpha,ß) = [11\text{-}(\alpha+ß), 9\text{-}(\alpha+ß)],\ \alpha \leq ß+1,\ \alpha+ß \leq 9$

$e_2^5\ [K_{14}\ (\alpha_1,ß_1)\ ] = [(\alpha_1+ß_1)+1,\ (\alpha_1+ß_1)\text{-}1] = K_{14}\ (\alpha_2,ß_2) = [11\text{-}(\alpha_2+ß_2),\ 9\text{-}(\alpha_2+ß_2)]$

$(\alpha_1+ß_1) + (\alpha_2+ß_2) = 10$

$K_{14}$ 's operation conditions are very strict $\alpha=ß$ or $\alpha=ß+1$ with $\alpha+ß \leq 9$, which generates two possible equivalences $R_3$

If $\alpha_1=ß_1$ and $\alpha_2=ß_2 \rightarrow (\alpha,\alpha)\ R_3^{11}\ (5\text{-}\alpha, 5\text{-}\alpha)$ $\quad 1 \leq \alpha \leq 4$

If $\alpha_1=ß_1+1$ and $\alpha_2=ß_2+1 \rightarrow (\alpha,\alpha\text{-}1)\ R_3^{12}\ (6\text{-}\alpha, 5\text{-}\alpha)$ $\quad 1 \leq \alpha \leq 5$

Crossed situations such as $\alpha_1 = ß_1$ and $\alpha_2 = ß_2+1$ are impossible

**Equivalences $R_3^{11}$**: $e_3^{11}\ (\alpha,\alpha) = (5\text{-}\alpha, 5\text{-}\alpha)\ \ 1 \leq \alpha \leq 4$

Thus, $(4,4)\ R_3^{11}\ (1,1)$ since $K_{14}\ (4,4) = (3,1)$ and $K_{14}\ (1,1) = (9,7)$, $(3,1)\ R_2^5\ (9,7)$

Associated functions: $e_2^2 \times e_3^{11} = e_3^{111}\ (\alpha,\alpha) = (5+\alpha, 5\text{-}\alpha)\ 1 \leq \alpha \leq 4$ (result from approach E)

The other products of functions $R_2$ with $R_3^{11}$ are impossible.

These equivalences allow to group several classes $C_2$ in classes $C_3$. Thus, for example,

$(4,4) \in C_2^5,\ (1,1) \in C_2^6,\ (4,4)\ R_3^{11}\ (1,1) \rightarrow C_2^5\ R_3\ C_2^6$, according to [13]

Similarly,

(4,4) $R_3^{11}$ (9,1) , (9,1) ∈ $C_2^6$ → $C_2^5$ $R_3$ $C_2^6$

(3,3) $R_3^{111}$ (2,2), (3,3) ∈ $C_2^{17}$, (2,2) ∈ $C_2^{16}$ → $C_2^{17}$ $R_3$ $C_2^{16}$

**Equivalences $R_3^{12}$**: $e_3^{12}$ (α, α-1) = (6-α, 5-α) 1≤α≤5

(4,3) $R_3^{12}$ (6-4, 5-4) = (2,1), $K_{14}$ (4,3) = (4,2), $K_{14}$ (2,1) = (8,6), (4,2) $R_2^5$ (8,6)

Associated functions:

$e_2^2$ x $e_3^{12}$ = $e_3^{121}$ (α,α-1) = (4+α, 5-α) 1 ≤ α ≤ 4

$e_2^7$ x $e_3^{12}$ = $e_3^{122}$ (α,α-1) = (5+α, 6-α) 1 ≤ α ≤ 4

Examples of class groupings:

(4,3) $R_3^{12}$ (2,1), (4,3) ∈ $C_2^4$, (2,1) ∈ $C_2^3$ → $C_2^4$ $R_3$ $C_2^3$

**Equivalences $R_3^2$** $e_3^2$ (α, t-α) = (t+5-α, α-5) 2α ≥t+1, t = 9, 10, 11

Equivalences $R_2^5$ are also generated by functions $K_5$, $K_6$, $K_9$ and $K_{17}$

$K_{17}$ (Table 2 in the Annex) and $R_2^5$ yield

$K_{17}$ (α,ß) = [(α-ß)+1, (α-ß)-1] α≥ß+1 9≤α+ß≤11; $e_2^5$ (α,ß) = (10-ß, 10-α)

α+ß = t, t=9, 10, 11 α-ß = 2α-t $K_{17}$ (α, t-α) = (2α-t+1, 2α-t-1), 2α-t ≥ 1

$e_2^5$ x $K_{17}$ ($\alpha_1$, t-$\alpha_1$) = (11+t-2$\alpha_1$, 9+t-2$\alpha_1$) = $K_{17}$ ($\alpha_2$, t-$\alpha_2$)

$\alpha_2$ = t+5-$\alpha_1$, $\alpha_1$+$ß_1$ = $\alpha_2$+$ß_2$=t, 2$\alpha_1$≥t+1, t = 9,10,11

**Equivalences $R_3^{21}$**: $e_3^{21}$ (α, 9-α) = (14-α, α-5), 5 < α < 9

They are the result of t=9

Example: (8,1) $R_3^{21}$ (6,3) →$C_2^3$ $R_3$ $C_2^4$, classes already related

Associated functions: $e_3^{211}$ = $e_2^2$ x $e_3^{21}$ (α, 9-α) = (α-4, α-5) 5 < α < 9

Examples: (8,1) $R_3^{211}$ (4,3), (7,2) $R_3^{211}$ (3,2), (6,3) $R_3^{211}$ (2,1)

Note that (7,2) $R_2^2$ (3,2) and, therefore, not interesting. The other two cases also group classes $C_2$ already related.

$e_3^{212}$ = $e_2^5$ x $e_3^{21}$ (α, 9-α) = (15-α, α-4) 5 < α < 9

Examples: (8,1) $R_3^{212}$ (7,4), (7,2) $R_3^{212}$ (8,3), (6,3) $R_3^{212}$ (9,2)

$e_3^{213}$ = $e_2^7$ x $e_3^{21}$ (α, 9-α) = (15-α, 14-α) 6 < α < 9

**Equivalences $R_3^{22}$**: $e_3^{22}$ (α, 10-α) = (15-α, α-5), 6 < α ≤ 9

They are the result of t = 10

Examples (9,1) $R_3^{22}$ (6,4), (8,2) $R_3^{22}$ (7,3) which relate already mentioned classes $C_2$

associated functions

The products are compatible with $e_2^2$ and $e_2^3$ in Group I and the corresponding ones in Group II, which coincide because α+ß =10. The product $e_2^5$ x $e_3^{22}$ = $e_3^{22}$. Naturally, there are no relations with new classes, since the product yields equivalent images $R_2$ to those yielded by $e_3^{22}$.

**Equivalences $R_3^{23}$**: $e_3^{23}$ (α, 11-α) = (16-α, α-5), 6 < α ≤ 9

They are the result of t=11. They do not provide new relations between classes.

Functions $K_1$, $K_2$, $K_5$ and $K_6$ generate equivalences $R_2^3$, while functions $K_5$ and $K_6$ also generate equivalences $R_2^2$.

**b) <u>$K_i \neq K_j$</u>**

Equivalences $R_3$ can be generated through equivalences $R_2$ between images of classes transformed by different functions $K_i$.

**Equivalences $R_3^3$**: $e_3^3$ (α, 4-α) = (α +(t+1)/2, (t-1)/2-α), t = 9, 11, 0 < α ≤ 4

They are the result of $e_2^3$ x $K_6$ ($α_1$,$ß_1$) = $K_{17}$ ($α_2$,$ß_2$)

**Equivalences $R_3^{31}$**: $e_3^{31}$ (α, 4-α) = (α +5, 4-α) 3 ≤ α ≤ 4

They are the result of t = 9

(3,1) $R_3^{31}$ (8,1) → $C_2^2$ $R_3$ $C_2^3$; (4,0) $R_3^{31}$ (9,0) → $C_2^8$ $R_3$ $C_2^7$

**Equivalences $R_3^{32}$**: $e_3^{32}$ (α, 4-α) = (α +6, 5-α) α = 3

They are the result of t = 11

(3,1) $R_3^{32}$ (9,2) → $C_2^2$ $R_3$ $C_2^3$

**Equivalences $R_3^4$**

Another example of these equivalences results from

$K_{17}$ ($α_1$,$ß_1$), $K_{14}$ ($α_2$,$ß_2$) and ($α_1$',$ß_1$') $R_2^5$ ($α_2$',$ß_2$')

Let us recall

$K_{17}$: $α_1 ≥ ß_1+1$, t= $α_1+ß_1$ $α_1-ß_1 = 2α_1-t$, t = 9, 10,11

$K_{14}$: $ß_2 ≤ α_2 ≤ ß_2+1$, $α_2+ß_2 ≤ 9$

Result: $(\alpha_1,ß_1)\ R_3^4\ (\alpha_2,ß_2) \rightarrow \alpha_1-ß_1 = \alpha_2+ß_2$

**Equivalence $\mathbf{R_3^{41}}$:** $e_3^{41}(\alpha, 10-\alpha) = (\alpha-5, \alpha-5)\quad \alpha>5$

$\alpha_2=ß_2,\quad 2\alpha_2 = 2\alpha_1 - t \rightarrow t = \dot{2} \rightarrow t = 10$

It generates the following relations

$(9,1)\ R_3^{41}\ (4,4) \rightarrow C_2^6\ \ R_3\ \ C_2^5$; $(8,2)\ R_3^{41}\ (3,3) \rightarrow C_2^{16}\ \ R_3\ \ C_2^{17}$

$(7,3)\ R_3^{41}\ (2,2) \rightarrow C_2^{17}\ \ R_3\ \ C_2^{16}$; $(6,4)\ R_3^{41}\ (1,1) \rightarrow C_2^5\ \ R_3\ \ C_2^6$

**Equivalence $\mathbf{R_3^{42}}$:** $e_3^{42}(\alpha, t-\alpha) = (\alpha-(t-1)/2,\ \alpha-(t+1)/2)\quad (t+1)/2 < \alpha$

$\alpha_2 = ß_2+1,\quad \alpha_2+ß_2 = 2\alpha_2-1 = \alpha_1-ß_1 = 2\alpha_1-t$

$\alpha_2 = \alpha_1 - (t-1)/2,\ ß_2 = \alpha_1 - (t+1)/2,\ t = 9,11$

**Equivalence $\mathbf{R_3^{421}}$:** $e_3^{421}(\alpha, 9-\alpha) = (\alpha-4,\alpha-5)\quad 4 < \alpha < 9$

It results from t = 9

$(8,1)\ R_3^{421}\ (4,3) \rightarrow C_2^3\ R_3\ \ C_2^4$; $(7,2)\ R_3^{421}\ (3,2) \rightarrow C_2^9\ \ R_3\ \ C_2^9$ (old)

$(6,3)\ R_3^{421}\ (2,1) \rightarrow C_2^4\ R_3\ \ C_2^3$; $(5,4)\ R_3^{421}\ (1,0) \rightarrow C_2^{12}\ \ R_3\ \ C_2^{11}$

Note that $(7,2)\ R_2^2\ (3,2)$

**Equivalence $\mathbf{R_3^{422}}$:** $e_3^{422}(\alpha, 11-\alpha) = (\alpha-5,\alpha-6)\quad 5 < \alpha \leq 9$

It results from t = 11

$(9,2)\ R_3^{422}\ (4,3) \rightarrow C_2^3\ R_3\ \ C_2^4$; $(8,3)\ R_3^{422}\ (3,2) \rightarrow C_2^9\ \ R_3\ \ C_2^9$ (old)

$(7,4)\ R_3^{422}\ (2,1) \rightarrow C_2^4\ R_3\ \ C_2^3$; $(6,5)\ R_3^{422}\ (1,0) \rightarrow C_2^{12}\ \ R_3\ \ C_2^{11}$

Table 6. Some equivalence functions $R_3$ in $A_4$

| Function $e_3$ | Domains of existence |
|---|---|
| $e_3^{11}$ (α,α) = (5-α, 5-α) | 1≤α≤4 |
| $e_3^{12}$ (α,α-1) = (6-α, 5-α) | 1≤α≤5 |
| $e_3^{21}$ (α,9-α) = (14-α, α-5) | 5<α<9 |
| $e_3^{22}$ (α,10-α) = (15-α, α-5) | 6<α≤9 |
| $e_3^{23}$ (α,11-α) = (16-α, α-5) | 6<α≤9 |
| $e_3^{31}$ (α,4-α) = (α+5, 4-α) | 3≤α≤4 |
| $e_3^{32}$ (α,4-α) = (α+6, 5-α) | α=3 |
| $e_3^{41}$ (α,10-α) = (α-5, α-5) | α>5 |
| $e_3^{42}$ (α,t-α) = (α-(t-1)/2, α-(t+1)/2) | α >(t-1)/2, t=9,11 |

In Table 6, some of the explained basic functions $e_3$ are summarized. Those resulting from the product with functions $e_2$ are not shown in the table, though some examples such as $e_3^{211} = e_2^{2} \times e_3^{21}$ or $e_3^{122} = e_2^{7} \times e_3^{2}$ have been given throughout the paper.

The product of functions $e_3$, given the strong restrictions of their arguments, are either inconsistent or only appliable to a few classes.

The explained equivalences $e_3$ are enough to group all classes $C_2$ in classes $C_3$ (Table 7). Any of the numbers m and n included in these classes will return the same third image.

$m \in C_3^i$ and $n \in C_3^i \rightarrow K^3(m) = K^3(n)$

Table 7. Equivalence classes $C_3$, $C_4$, $C_5$ and $C_6$ in $A_4$

| $R_3$ | | | $R_4$ | | | $R_5$ | | | $R_6$ | | |
|---|---|---|---|---|---|---|---|---|---|---|---|
| Classes $C_3$ | Classes $C_2$ included | Image of $C_3$ | Classes $C_4$ | Classes $C_3$ included | Image of $C_4$ | Classes $C_5$ | Classes $C_4$ included | Image of $C_5$ | Classes $C_6$ | Classes $C_5$ included | Image of $C_6$ |
| $C_3^1$ | $C_2^1, C_2^2$<br>$C_2^3, C_2^4$ | (6,2) | $C_4^1$ | $C_3^1, C_3^2$<br>$C_3^3, C_3^4$ | (6,2) | $C_5^1$ | $C_4^1, C_4^2$<br>$C_4^3, C_4^4$ | (6,2) | $C_6^1$ | $C_5^1, C_5^2$<br>$C_5^3, C_5^4$ | (6,2) |
| $C_3^2$ | $C_2^5, C_2^6$ | $C_2^2$ | | | | | | | | | |
| $C_3^3$ | $C_2^7$ | $C_2^3$ | | | | | | | | | |
| $C_3^4$ | $C_2^8$ | $C_2^4$ | | | | | | | | | |
| $C_3^5$ | $C_2^9$ | $C_2^5$ | $C_4^2$ | $C_3^5, C_3^6$ | $C_3^2$ | | | | | | |
| $C_3^6$ | $C_2^{10}$ | $C_2^6$ | | | | | | | | | |
| $C_3^7$ | $C_2^{11}, C_2^{12}$ | $C_2^7$ | $C_4^3$ | $C_3^7$ | $C_3^3$ | | | | | | |
| $C_3^8$ | $C_2^{13}$ | $C_2^8$ | $C_4^4$ | $C_3^8$ | $C_3^4$ | | | | | | |
| $C_3^9$ | $C_2^{14}$ | $C_2^9$ | $C_4^5$ | $C_3^9$ | $C_3^5$ | $C_5^2$ | $C_4^5$ | $C_4^2$ | | | |
| $C_3^{10}$ | $C_2^{15}$ | $C_2^{12}$ | $C_4^6$ | $C_3^{10}$ | $C_3^7$ | $C_5^3$ | $C_4^6$ | $C_4^3$ | | | |
| $C_3^{11}$ | $C_2^{16}, C_2^{17}$ | $C_2^{13}$ | $C_4^7$ | $C_3^{11}$ | $C_3^8$ | $C_5^4$ | $C_4^7$ | $C_4^4$ | | | |
| $C_3^{12}$ | $C_2^{18}$ | $C_2^{14}$ | $C_4^8$ | $C_3^{12}$ | $C_3^9$ | $C_5^5$ | $C_4^8$ | $C_4^5$ | $C_6^2$ | $C_5^5$ | $C_5^2$ |
| $C_3^{13}$ | $C_2^{19}$ | $C_2^{15}$ | $C_4^9$ | $C_3^{13}$ | $C_3^{10}$ | $C_5^6$ | $C_4^9$ | $C_4^6$ | $C_6^3$ | $C_5^6$ | $C_5^3$ |
| $C_3^{14}$ | $C_2^{20}$ | $C_2^{16}$ | $C_4^{10}$ | $C_3^{14}$ | $C_3^{11}$ | $C_5^7$ | $C_4^{10}$ | $C_4^7$ | $C_6^4$ | $C_5^7$ | $C_5^4$ |

On the approach used

a) New equivalences

It is interesting to find new equivalences ([5] and [5a]), but approach used (approach C) generates old equivalences ([6]) as well. Thus, (4,4) $R_3^{11}$ (1,1) is a new equivalence because (4,4) $\not R_2$ (1,1)

4004 (4,4) → 4356 (3,1) → 3087 (8,4) → 8352

1100 (1,1) → 1089 (9,7) → 9621 (8,4) → 8352

But (7,2) $R_3^{211}$ (3,2) is old because (7,2) $R_2^2$ (3,2)

2953 (7,2) → 7173 (6,4) → 6354

2425 (3,2) → 3177 (6,4) → 6354

$K^3$ (2953) = $K^3$ (2425) but also $K^2$ (2953) = $K^2$ (2425)

These old equivalences occur for those classes (α,ß) where $e_3 = e_2$.

In our example

$e_3^{211} = e_2^2 \times e_3^{21} = e_2^2 \rightarrow e_3^{21}$ (α,ß) = (α,ß)

Which actually does occur

$e_3^{21}$ (α,9-α) = (14-α, α-5) → $e_3^{21}$ (7,2) = (7,2)

Old equivalences $R_1$ can also occur. Thus,

$K_i$ ($α_1$,$ß_1$) = (α'$_1$,ß'$_1$), e (α'$_1$,ß'$_1$) = (α'$_1$,ß'$_1$) and $K_j$ ($α_2$,$ß_2$) = (α'$_1$,ß'$_1$) → ($α_1$,$ß_1$) = ($α_2$,$ß_2$) o ($α_1$,$ß_1$) $R_2$ ($α_2$,$ß_2$)

This happens when equivalence $e_2^5$ and functions $K_i$ and $K_j$ capable of generating class (6,4) are used, since $e_2^5$ (6,4) = (6,4)

Such functions can have the same inverse image, like $K_6$ and $K_{14}$, or inverse images $R_2$, like $K_{14}$ and $K_5$

$K_6$ (3,2) = (6,4), $K_{14}$ (3,2) = (6,4), $K_5$ (7,2) = (6,4), (3,2) $R_2^2$ (7,2)

b) Wrong equivalences

- If $K_i$ ($α_1$,$ß_1$) = (α'$_1$, 0) and $K_j$ ($α_2$,$ß_2$) = (α'$_2$,ß'$_2$), ß'$_2$ > 0
  then applying an equivalence $R_2$ is absurd, as classes (α,0) are only equivalent to other classes (11-α, 0).
  This situation can be addressed by restraining the equivalence's domain of existence, removing the prescribed cases. Thus, $R_3^{21}$ at first limited to $5 \leq α \leq 9$, leads to (9,0) $R_3^{21}$ (5,4) which is wrong since K (5,4) = (2,0) is inconsistent with K (9,0) = (8,1).The solution chosen is $R_3^{21}$ $5 < α < 9$, since $R_3^{21}$ is symmetric (5,4) $R_3^{21}$ (9,0), thus removing the inconsistent situations. The wrong specific result (9,0) $R_3^{21}$(54) occurs because $R_3^{21}$ uses $K_i = K_j = K_{17}$ and $e_2^5$. $K_{17}$ cannot operate on (9,0), yielding the wrong result $K_{17}$ (9,0) = (10,8). But $e_2^5$ (10,8) = (2,0) which coincides with the correct image $K_{17}$ (5,4) = (2,0). Schematically,

$$(9,0) \xrightarrow{K_{17}} (10,8) \xrightarrow{e_2^5} (2,0)$$

Error of domain    Error of domain

$$(5,4) \xrightarrow{K_{17}} (2,0)$$

⇒ (9,0) $\not{R}_3^{21}$(5,4)

- If K ($α_1$,$ß_1$) = (α'$_1$, 0) and K ($α_2$,$ß_2$) = (α'$_2$, 0) then
  (α'$_1$, 0) $R_2^1$ (α'$_2$, 0) → ($α_1$,$ß_1$) $R_3$ ($α_2$,$ß_2$)
  The only valid equivalence $R_2$ that can be applied in this case is $e_2^1$ (α,0) = (11-α,0). However, there is no $R_3$ of this type. There are though some $R_4$, as will be shown later.

## 7.4 Equivalence $R_4$ and above

As previously mentioned, as the equivalence order is increased, the **new equivalence** relations only occur between a very reduced number of classes. The most striking example analyzed in the previous section is $R_3^{32}$, which only relates two classes (α,ß).

**Approach C**

**Equivalences $R_4^1$**: $e_4^1$ (α,0) = (α-5,0), 6 < α ≤ 9

It is based on $K_{25}$, $K_{26}$ and $e_3^{21}$ → $\alpha_1$-$\alpha_2$ = 5

(9,0) $R_4^1$ (4,0) → $C_2^7$ $R_4$ $C_2^8$ → $C_3^4$ $R_4$ $C_3^4$

(8,0) $R_4^1$ (3,0) actually is (8,0) $R_2^1$ (3,0)

(7,0) $R_4^1$ (2,0) → $C_2^8$ $R_4$ $C_2^7$ → $C_3^8$ $R_4$ $C_3^7$

Associated functions:

$e_4^{11}$ = $e_2^1$ x $e_4^1$ (α,0) = (16-α, 0) , 6 < α ≤ 9

(9,0) $R_4^{11}$ (7,0), (7,0) $R_4^{11}$ (9,0) → $C_2^8$ $R_4$ $C_2^7$ → $C_3^4$ $R_4$ $C_3^4$

Applying $e_2^1$ to the argument of $R_4^1$ yields

$e_4^{12}$ (11-α,0) = (α-5, 0) , 6 < α ≤ 9

(2,0) $R_4^{12}$ (4,0), (3,0) $R_4^{12}$ (3,0), (4,0) $R_4^{12}$ (2,0)

**Equivalences $R_4^2$**: $e_4^2$ (α, t-α) = (5,5), α = (t+5)/2, t = 9,11

It is based on $K_{17}$, $K_{10}$ and $e_3^{41}$ → $\alpha_2$ = $ß_2$ = 5

**Equivalences** $R_4^{21}$: $e_4^{21}$ (α, 9-α) (5,5), α = 7

It results from t=9

(7,2) $R_4^{21}$ (5,5) → $C_2^9$ $R_4$ $C_2^{10}$ → $C_3^5$ $R_4$ $C_3^6$

**Equivalences** $R_4^{22}$: $e_4^{22}$ (α, 11-α) = (5,5), α = 8

It results from t=11

(8,3) $R_4^{22}$ (5,5) → $C_2^9$ $R_4$ $C_2^{10}$ → $C_3^5$ $R_4$ $C_3^6$

(7,2) $R_4^{221}$ (5,5) results from (8,3) $R_2^5$ (5,5) and [9]

Similarly with (3,2) and (8,7)

**Equivalences $R_4^3$**: $e_4^3$ (6,4) = (2,0)

It uses $K_{17}$, $K_{26}$ and $e_3^{31}$ → t=10, $\alpha_1$ = (t+2)/2 = 6, $\alpha_2$ = 4-2$\alpha_1$+t = 2

$e_4^3$ ($\alpha_1$, t-$\alpha_1$) = ($\alpha_2$, 0)

The equivalence is limited to only one case

(6,4) $R_4^3$ (2,0) → $C_2^5$ $R_4$ $C_2^7$ → $C_3^3$ $R_4$ $C_3^3$

The secondary equivalences result from replacing (6,4) and/or (2,0) by equivalences $R_2$ o $R_3$, under [9]. Thus, using $R_2$

(4,4) $R_4^{31}$ (2,0), (4,4) $R_4^{32}$ (9,0), (6,6) $R_4^{33}$ (2,0), (6,6) $R_4^{34}$ (9,0)

And similarly with $R_3$. (4,4) $R_3^{11}$ (1,1) yields (1,1) $R_4^{35}$ (2,0), etc.

**Equivalence $\mathbf{R_5^1}$**: (1,0) $R_5^1$ (7,5)

It results from $K_{26}$, $K_1$ and $e_4^1$

$e_4^1$ [$K_{26}$ ($\alpha_1$,0)] = (5-$\alpha_1$, 0), $\alpha_1$-1 = 0

$K_1$ ($\alpha_2$,$ß_2$) = (2$\alpha_2$-10, 2$ß_2$-10) 6≤ $\alpha_2$≤9

If made equal it returns 2 $\alpha_2$ + $\alpha_1$ = 15, 2$ß_2$-10 = 0 $\alpha_1$ = 1

There is only one solution (1,0) $R_5^1$ (7,5) → $C_2^{11}$ $R_5$ $C_2^{13}$ → $C_3^7$ $R_5$ $C_3^8$

And all the associated functions by replacing (1,0) and/or (7,5) by equivalents of a lower range

**Approach D**

**Equivalence $\mathbf{R_5^1}$**: (1,0) $R_5^1$ (7,5)

At this stage, where some equivalences of a high range operate only on one pair, it sounds logical to use approach D.

For example, applying it to (1,0) $R_5^1$ (7,5) is impossible since (1,0) does not have an inverse image

But it can be applied to (9,0) $R_4^1$ (4,0)

$K_{26}^{-1}$ (9,0) = (1,0), $K_2^{-1}$ (4,0) = (7,5) → (1,0) $R_5^2$ (7,5)

**Equivalence $\mathbf{R_6^1}$**: (8,2) $R_6^1$ (6,0)

The non-existence of an inverse image may be ignored if a given class has an equivalence of a lower order with another class

For instance, the existence of (9,0) $R_4^1$ (4,0) has enabled us to find (1,0) $R_5^2$ (7,5) and, as previously mentioned, approach D cannot be iterated on classes (1,0) and (7,5). But because (1,0) $R_3^4$ (5,4) then, under [9] (7,5) $R_5$ (5,4), $K^{-1}$ (7,5) = $K_{17}^{-1}$ (7,5)= = (8,2), $K_{25}^{-1}$ (5,4) = (6,0) from where

(8,2) $R_6^1$ (6,0) → $C_2^{16}$ $R_6$ $C_2^{15}$ → $C_3^{11}$ $R_6$ $C_3^{10}$

**Equivalence $R_7^1$**: (9,4) $R_7^1$ (8,5)

$K_2^{-1}$ (8,2) = (9,4), $K_2^{-1}$ (6,0) = (8,5) as (8,2) $R_6^1$ (6,0)

(9,4) $R_7^1$ (8,5) → $C_2^{20}$ $R_7$ $C_2^{19}$ → $C_3^{14}$ $R_7$ $C_3^{13}$

The process cannot be iterated, (9,4) and (8,5) do not have inverse images. We have reached the top of transformations. All the transformed classes of (9,4) and (8,5) with all their equivalences of order 6 and lower, converge at (6,2).

**<u>Sequential approach F</u>**

**Equivalences $R_5^2$, $R_6^2$ and $R_7^2$**

Let us see an example of operators $K^r$ being used sequentially (Approach F)

Let the process be

$$(\alpha,\beta) \xrightarrow{K_1} (\alpha', \beta') \xrightarrow{K_{25}} (\alpha'', \beta'') \xrightarrow{K_{17}} (\alpha''', \beta''')$$

$K_1$ (α, t-α) = (α', ß') = (2α-10, 2t-2α-10) 6≤α≤9, 11≤t≤17

The compatibility of $K_{25}$ x $K_1$ = $K^2$

$K_{25}$ (α,0) = (α-1, 10-α), 6 ≤α≤ 9 demands 2t-2α-10 = 0 ↔ α = t-5

That is,

$K_1$(t-5,5) = (2t-20, 0) 11 ≤ t ≤ 14

$K^2$ = $K_{25}$ x K1 (t-5, 5) = (2t-21, 30-2t) t = 13,14

$K_{17}$ (α, $t_1$-α) = (2α-$t_1$+1, 2α-$t_1$-1), 2α-$t_1$ ≥ 1

Its compatibility with $K_2$ demands 30-2t = $t_1$ – (2t-21) ↔ $t_1$ = 9

Therefore

$K_{17}$ (α,9-α) = (2α-8, 2α-10)

$K^3$ = $K_{17}$ x $K^2$ (t-5, 5) = (4t-50, 4t-52) t= 13,14

The process is thus

$$(t\text{-}5, 5) \xrightarrow{K=K_1} (2t\text{-}20, 0) \xrightarrow{K^2=K_{25}xK_1} (2t\text{-}21, 30\text{-}2t) \xrightarrow{K^3=K_{17}xK_{25}xK_1} (4t\text{-}50, 4t\text{-}52) \quad t=13,14$$

The process can split in two, the previous one for t = 13 and, if t is replaced with t+1

$$\begin{array}{cccc} & K & K^2 & K^3 \\ (t\text{-}4, 5) & \rightarrow (2t\text{-}18, 0) & \rightarrow (2t\text{-}19, 28\text{-}2t) & \rightarrow (4t\text{-}46, 4t\text{-}48) \end{array}$$

also for t=13

If we prove that (4t-46, 4t-48) $R_r$ (4t-50, 4t-52)

then necessarily

(2t-19, 28-2t) $R_{r+1}$ (2t-21, 30-2t), (2t-18, 0) $R_{r+2}$ (2t-20, 0), (t-4, 5) $R_{r+3}$ (t-5, 5)

with t = 13

Defining values t=13

Process 1: $K^3$ (8,5) = (2,0), $K^2$ (8,5) = (5,4), K (8,5) = (6,0), $K^0$ (8,5) = (8,5)

Process 2: $K^3$ (9,5) = (6,4), $K^2$ (9,5) = (7,2), K (9,5) = (8,0), $K^0$ (9,5) = (9,5)

Since we have proven that (6,4) $R_4^2$ (2,0) then there is

**Equivalences $\mathbf{R_5^2}$** : (2t-19, 28-2t) $R_5$ (2t-21, 30-2t) $\xrightarrow{t=13}$ (7,2) $R_5$ (5,4)

**Equivalences $\mathbf{R_6^2}$**: (2t-18, 0) $R_6$ (2t-20, 0) → (8,0) $R_6$ (6,0)

**Equivalences $\mathbf{R_7^2}$**: (t-4, 5) $R_7$ (t-5, 5) → (9,5) $R_7$ (8,5)

Classes (9,5) and (8,5) do not have an inverse image. All classes resulting from their transformations will verify $K^6(\alpha,\beta) = (6,2)$ and the numerical image K(6,2) = 6174

We have then proven that equivalences $R_7$ which relate these classes without an inverse image exist. As it was also proven for the rest of classes without an inverse image such as (9,4), (1,0) or (5,5), we can thus state that equivalences $R_7$ encompass **all** classes and therefore the system converges towards 6174.

Table 7 shows equivalence classes $C_4$ and above.

Graph 4 shows transformations of the higher order classes.

## 8. Equivalences in $A_3$

By applying the notions of equivalence to 3-digit numbers, we can know how equivalence classes, sets of numbers with shared images, start to group and create the branches and leaves in the transformation tree.

Functions in $A_3$ are:

$f(\alpha) = (\alpha\text{-}1\ \ 9\ \ 10\text{-}\alpha)\ \ 0 < \alpha \leq 9$

$K_1(\alpha) = (\alpha\text{-}1),\ \ 6 \leq \alpha \leq 9$

$K_2(\alpha) = (10\text{-}\alpha),\ \ 0 < \alpha \leq 5$

and balance

$(\alpha_E) = (5),\ \ f(5) = n_E = 495$

**Equivalences $\mathbf{R_0}$ and $\mathbf{R_1}$**

- $C_0^i$: each class contains only one number
- $C_1^i$: each class consists of all numbers yielding the same image and for what is stated in section 5b, $C_1 = (\alpha)$.

**Equivalences $\mathbf{R_2}$**

- $C_2^1$: they consist of all numbers m and n – and corresponding classes ($\alpha_1$) and ($\alpha_2$) – such that
  m $R_2$ n ↔ m' $R_1$ n' ↔ ($\alpha'_1$) = ($\alpha'_2$)
  thus including the identity
  $e_2^0$ ($\alpha$) = $\alpha$, de $K_i$ ($\alpha_1$) = $K_i$ ($\alpha_1$)
  or 11's supplementation when $K_1(\alpha_1) = K_2(\alpha_2)$
  $e_2^1$ ($\alpha$) = (11-$\alpha$) de $\alpha_1$-1 = 10-$\alpha_2$ $0<\alpha\leq 9$

  In the same way that happened with equivalences $e_2^0$ and $e_2^1$ in $A_4$, these two functions have an isomorphic abelian group algebraic structure to $\mathbb{Z}_2$ and the permutational symmetric group $S_2$.

  Equivalences $R_1$ and $R_2$ generate the following class graph

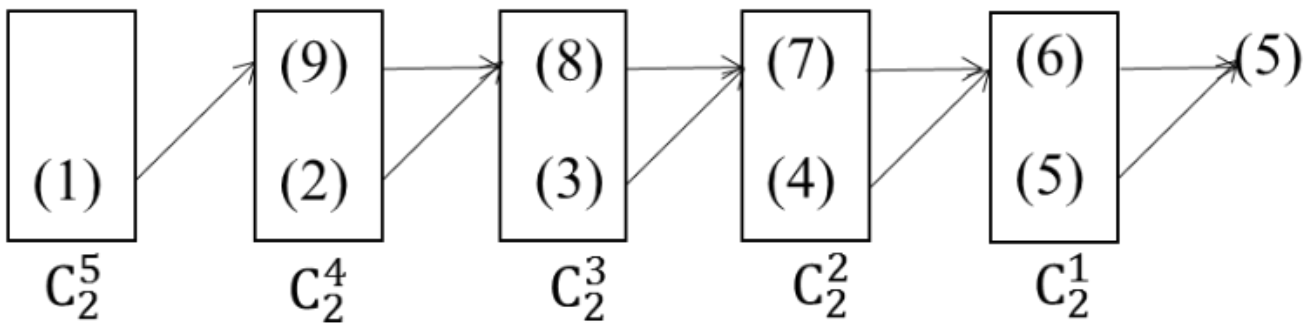

  - If $6 \leq \alpha \leq 9$, the distance D or number of transformations to reach the balance class (5) will be D = $\alpha$-5
  - If $1 \leq \alpha < 5$, D = 6-$\alpha$
- Equivalences $R_3$ and higher

This tree is not compatible with new **$R_3$**.

In each step towards a higher equivalence all the previous equivalences (old) last, and the class containing (5) merges with the adjoining one. For example (7) $R_3$ (6) since $K^2(7) = K^2(6) = (5)$

$C_3^1 = \{ C_2^1,\ C_2^2\},\ C_3^2 = C_2^3,\ C_3^3 = C_2^4,\ C_3^4 = C_2^5$

This way with $R_6$ there is only one equivalence class $C_6^1$ containing all previous classes.

## 9. <u>Equivalences in $A_2$</u>

In set $A_2$, which consists of 2-digit numbers, transformations are determined by a single parameter, just as it happens in $A_3$.

Functions in $A_2$ are

f($\alpha$) = ($\alpha$-1 10-$\alpha$), $0 < \alpha \leq 9$

$K_1(\alpha)$ = (2$\alpha$-11), $6 \leq \alpha \leq 9$

$K_2(\alpha)$ = (11-2$\alpha$), $0 < \alpha \leq 5$

Balance does not exist, but a 5-link cycle. The 9 classes $C_1$ = ($\alpha$) graph is

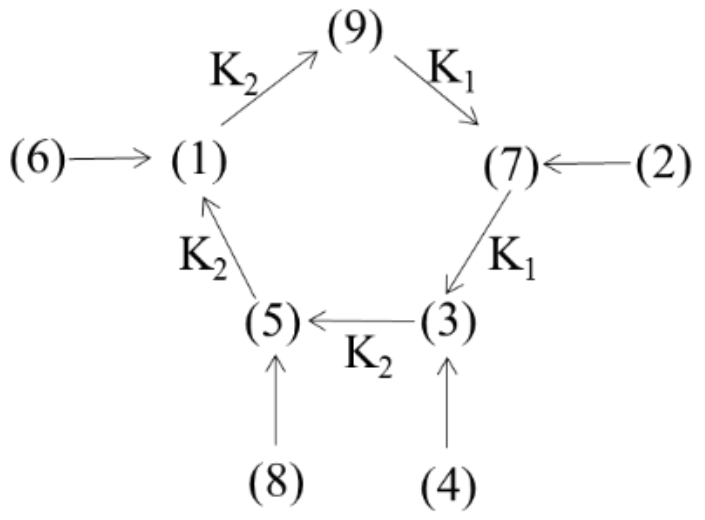

$R_2^0$: α'1-1 = α'2-1, 10-α'1 = 10-α'2 ↔ $e_2^0$ (α) = (α)

and 11's supplementation

$R_2^1$: α'1-1 = 10-α'2, 10-α'1 = α'2-1↔ $e_2^1$ (α) = (11-α), $2 \leq \alpha \leq 9$

These equivalences have the same isomorph group structure to $\mathbb{Z}_2$ previously mentioned for $A_4$ and $A_3$.

The equivalence classes $C_2$ are

$C_2^1$ = {(9), (2)}, $C_2^2$ ={(7),(4)}, $C_2^3$={(8),(3)}, $C_2^4$={(6),(5)}, $C_2^5$={(1)}

linked by graph

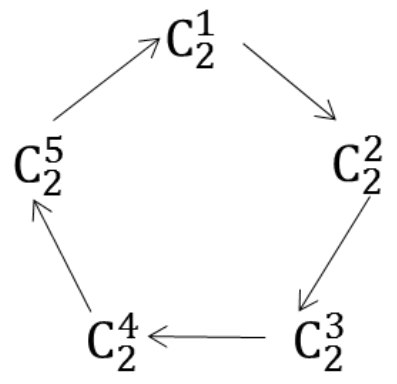

**10. Equivalences in $A_5$**

Set $A_5$, consisting of numbers of up to five digits, has the following basic functions

$f_1$ (α,ß) = (α ß-1 9 9-ß 10-α), α ≥ ß, 0 < ß ≤ 9

$f_2$ (α,0) = (α-1 9 9 9 10-α), 1 ≤ α ≤ 9

Functions $K_i$ are shown in the Annex (Table 4).

**Equivalences $R_0$**

m $R_0$ m ∀m∈ $A_5$

There are 99.990 classes $C_0$ each of them composed of just one number

**Equivalences $R_1$**

m $R_1$ n ↔ K(m) = K(n) ↔ p(m) = p(n) ∀m,n ∈ $A_5$

m ∈ ($\alpha_1$,$ß_1$) and m ∈ ($\alpha_2$,$ß_2$), m $R_1$ n ↔ ($\alpha_1$,$ß_1$) = ($\alpha_2$,$ß_2$)

There are 54 classes $C_1$, as in $A_4$.

Graph 5 shows the transformation tree corresponding to classes $C_1 = (\alpha,ß)$. There is one 2-link cycle and two 4-link cycles.

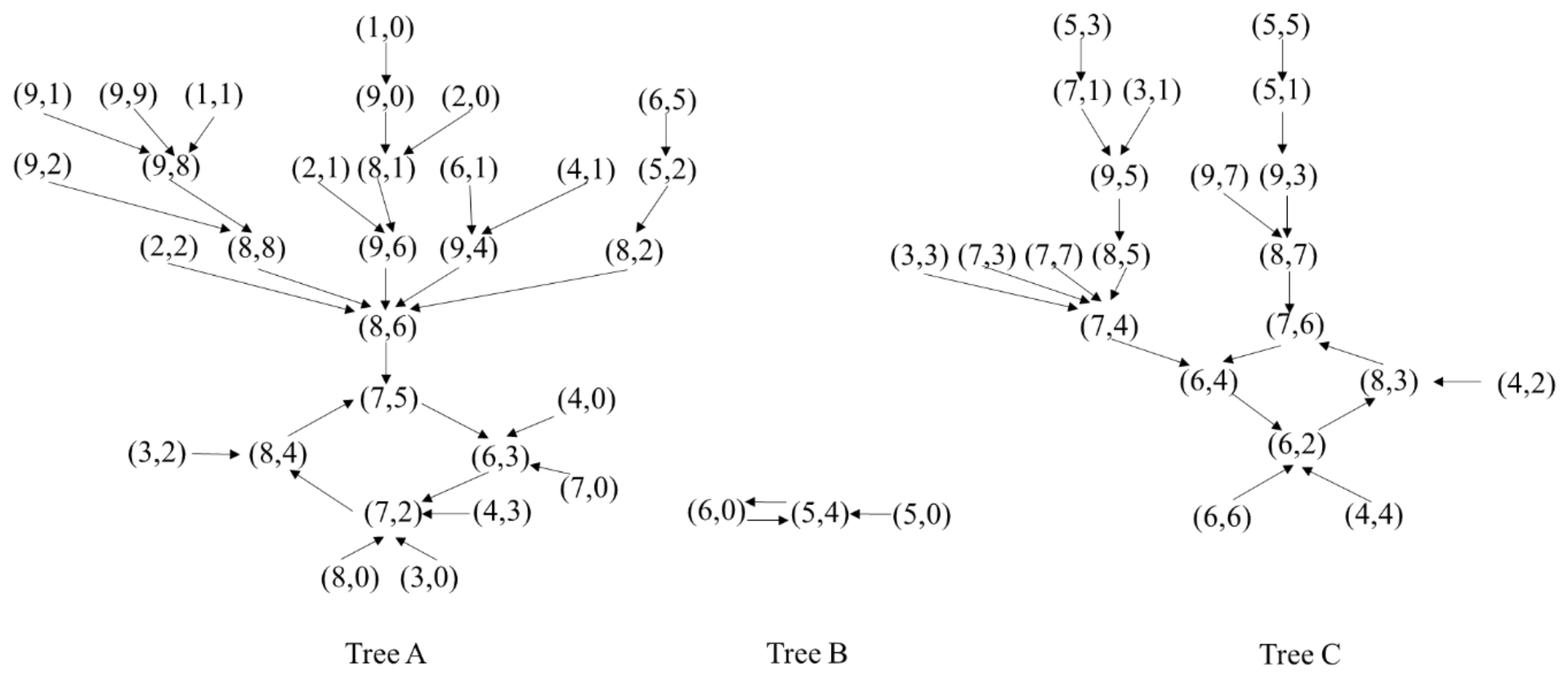

Graph 5. Parametric transformations in $A_5$

**Equivalences $R_2$**

$(\alpha_1,ß_1)$ $R_2$ $(\alpha_2,ß_2) \leftrightarrow (\alpha'_1,ß'_1) = (\alpha'_2,ß'_2)$

Equivalences $R_2^1$, $R_2^2$, $R_2^3$ and $R_2^4$, previously defined in $A_4$, remain

$\mathbf{R_2^1}$: $e_2^1$ $(\alpha,0) = (11\text{-}\alpha,0)$, derives from $K_{25}$ $(\alpha_1,0) = K_{26}$ $(\alpha_2,0)$

$\mathbf{R_2^2}$: $e_2^2$ $(\alpha,ß) = (10\text{-}\alpha,ß)$, $\alpha+ß \leq 10$, derives for example from $K_{14}$ $(\alpha_1,ß_1) = K_{17}$ $(\alpha_2,ß_2)$

$\mathbf{R_2^3}$: $e_2^3$ $(\alpha,ß) = (\alpha,10\text{-}ß)$, $\alpha+ß \geq 10$, derives for example from $K_9$ $(\alpha_1,ß_1) = K_{17}$ $(\alpha_2,ß_2)$

$\mathbf{R_2^4}$: $e_2^4$ $(\alpha,ß) = (10\text{-}\alpha,10\text{-}ß)$, $\alpha=ß$, derives from $K_9$ $(\alpha_1,ß_1) = K_{14}$ $(\alpha_2,ß_2)$

Equivalences $R_2^5$, $R_2^6$, $R_2^7$ and $R_2^8$ disappear, except when coinciding with the previous ones.

As in $A_4$, $A_3$ and $A_2$ equivalences $e_2^0$ and $e_2^1$ have an isomorphic group structure in $\mathbb{Z}_2$. As in $A_4$ equivalences $e_2^0$, $e_2^2$, $e_2^3$ and $e_2^4$ have an isomorphic group structure of Klein group.

New equivalences appear

$\mathbf{R_2^5}$: $e_2^5(\alpha,ß) = [1/2\ (\alpha+ß+1),\ 11\text{-}\alpha]$, $\alpha \geq ß+1$, $5 \leq ß \leq 8$, $9 \leq \alpha+ß \leq 11$

From $K_1(\alpha_1,ß_1)=K_{17}$ $(\alpha_2,ß_2)$

Equivalences $R_1$ and $R_2^3$ are generated, as well as new equivalences such as (8,5) $R_2^5$ (7,3) and (9,6) $R_2^5$ (8,2)

$\mathbf{R_2^6}$: $e_2^6$ (α, 11-2α) = (11-2α, 0) only valid for (4,3) $R_2^6$ (3,0)

From $K_{14}$ ($α_1$,$ß_1$)=$K_{26}$ ($α_2$,$ß_2$)

$\mathbf{R_2^7}$: $e_2^7$ (α, 19-2α) = (11-α, 0) only valid for (7,5) $R_2^7$ (4,0)

From $K_1$ ($α_1$,$ß_1$)=$K_{26}$ ($α_2$,0)

There are no more new different $R_2$, except for those resulting from the transitivity (product) between $R_2$. Thus,

(7,5) $R_2^7$ (4,0) and (4,0) $R_2^1$ (7,0) → (7,5) $R_2^{71}$ (7,0)

(9,6) $R_2^5$ (8,2) and (9,4) $R_2^3$ (9,6), (8,2) $R_2^2$ (2,2) → (9,4) $R_2^{51}$ (2,2)

Table 8 shows the partition in classes $C_2$ and Graph 6 shows the corresponding transformation trees.

Table 8. Equivalence classes $C_2$ in $A_5$, Tree A

| Equivalence classes $C_2$ | $C_2^1$ | $C_2^2$ | $C_2^3$ | $C_2^4$ | $C_2^5$ | $C_2^6$ | $C_2^7$ | $C_2^8$ | $C_2^9$ | $C_2^{10}$ | $C_2^{11}$ | $C_2^{12}$ | $C_2^{13}$ |
|---|---|---|---|---|---|---|---|---|---|---|---|---|---|
| Classes $C_1$ included | (3,0) (4,3) (6,3) (8,0) | (3,2) (7,2) | (4,0) (7,0) (7,5) | (8,4) (8,6) | (2,2) (8,2) (8,8) (9,4) (9,6) | (9,2) (9,8) | (2,1) (8,1) | (5,2) | (4,1) (6,1) | (1,1) (9,1) (9,9) | (2,0) (9,0) | (6,5) | (1,0) |
| Image of $C_2$ | (7,2) | (8,4) | (6,3) | (7,5) | (8,6) | (8,8) | (9,6) | (8,2) | (9,4) | (9,8) | (8,1) | (5,2) | (9,0) |

| Tree B | | Tree C | | | | | | | | | | |
|---|---|---|---|---|---|---|---|---|---|---|---|---|
| $C_2^{14}$ | $C_2^{15}$ | $C_2^{16}$ | $C_2^{17}$ | $C_2^{18}$ | $C_2^{19}$ | $C_2^{20}$ | $C_2^{21}$ | $C_2^{22}$ | $C_2^{23}$ | $C_2^{24}$ | $C_2^{25}$ | $C_2^{26}$ |
| (5,0) (6,0) | (5,4) | (4,2) (6,2) | (4,4) (6,4) (6,6) | (7,4) (7,6) | (3,3) (7,3) (7,7) (8,5) | (8,3) (8,7) | (9,5) | (9,3) (9,7) | (3,1) (7,1) | (5,1) | (5,3) | (5,5) |
| (5,4) | (6,0) | (8,3) | (6,2) | (6,4) | (7,4) | (7,6) | (8,5) | (8,7) | (9,5) | (9,3) | (7,1) | (5,1) |

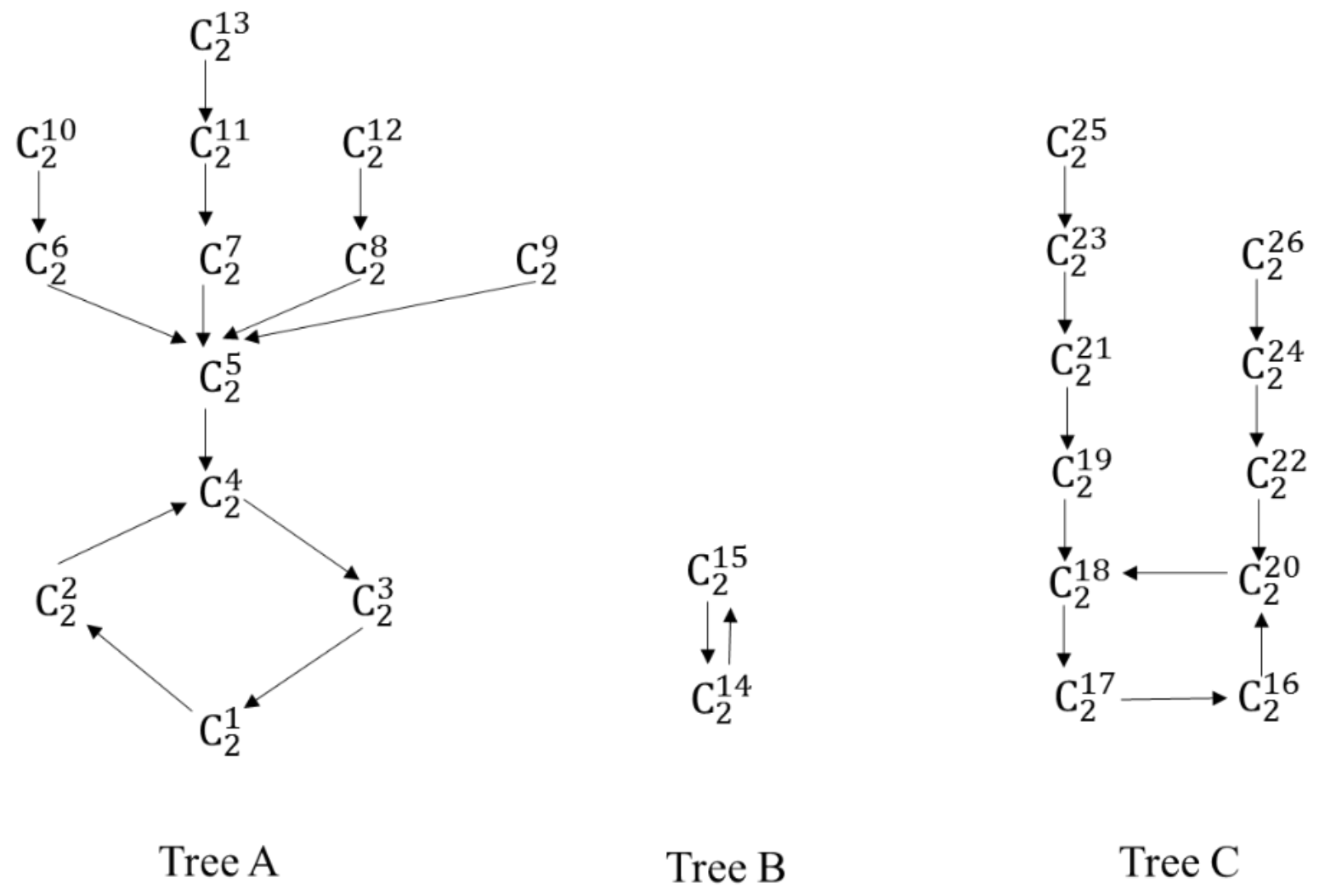

Graph 6. Flow of transformations between classes $C_2$ en $A_5$

**Equivalences $R_3$**

New equivalences need two links in the chains to get to the shared class α,ß. Let us use approach C

$\mathbf{R_3^1}$: $e_3^1$ (α,ß) = [1/2 (29 – (α+ß)), 11- α] valid for (9,6), (8,7) and (8,5)

Derives from $e_2^3$ x $K_1(\alpha_1,ß_1) = (\alpha_1-1, 19-\alpha_1-ß_1) = K_{17}$ $(\alpha_2-ß_2)$

(9,6) $R_3^1$ (7,2), (8,7) $R_3^1$ (7,3), (8,5) $R_3^1$ (8,3)

$\mathbf{R_3^2}$: $e_3^2$ (α,ß) = (α, 2α+ß-18), α ≥ 8, ß ≥ 5, α ≥ ß+1

Derives from $e_2^3$ x $K_1(\alpha_1,ß_1) = K_2$ $(\alpha_2,ß_2)$

(8,7) $R_3^2$ (8,5), (8,6) $R_2^3$ (8,4), (8,5) $R_3^2$ (8,3)

For α=9, $e_3^2$(α,ß) = (α,ß) → $e_3^2 = e_1 = e_2^0$

$\mathbf{R_3^3}$: $e_3^3$ (α,ß) = (2α-ß-11, 11-α), α = 9, 3 ≤ ß ≤ 5

Derives from $e_2^3$ x $K_2(\alpha_1,ß_1) = [\alpha_1-1, 9-(\alpha_1-ß_1)] = K_6$ $(\alpha_2,ß_2)$

(9,4) $R_3^3$ (3,2), (9,3) $R_3^3$ (4,2), (9,2) $R_3^3$ (5,2)

$\mathbf{R_3^4}$: $e_3^4$ (α,ß) = (5-α, ß) valid for α = 2, ß = 2

Derives from $e_2^3$ x $K_{14}(\alpha_1,ß_1) = (10-ß_1, 2\alpha_1) = K_{14}$ $(\alpha_2,ß_2)$

(2,2) $R_3^4$ (3,2)

$\mathbf{R_3^5}$: $e_3^5$ (α,ß) = (9-2α-ß, ß) valid for α = 2, ß = 2 valid for α = 2, ß = 1,2

Derives from $e_2^3$ x $K_{14}(\alpha_1,ß_1) = (10-ß_1, 2\alpha_1) = K_6$ $(\alpha_2,ß_2)$

(2,2) $R_3^5$ (3,2), (2,1) $R_3^5$ (4,1)

$R_3^5$ coincides with $R_3^4$ when α=4-ß

$R_3^5$ coincides with $R_3^1$ when α=ß+1

$\mathbf{R_3^6}$: $e_3^6$ (α,ß) = (8-α-2ß, ß) valid for α = 2, ß = 1,2

Derives from $e_2^3$ x $K_6(\alpha_1,ß_1)$ = $(10-ß_1, \alpha_1+ ß_1+1)$ = $K_6(\alpha_2,ß_2)$, (2,1) $R_3^6$ (4,1)

$\mathbf{R_3^7}$: $e_3^7$ (α,ß) = (8-α/2-2ß, α /2+ß) valid for α = 2, ß = 1

Derives from $e_2^5$ x $K_6(\alpha_1,ß_1)$ = $[1/2(20- \alpha_1-2ß_1), 1+ß_1]$ = $K_6(\alpha_2,ß_2)$

(2,1) $R_3^7$ (5,2)

$\mathbf{R_3^8}$: $e_3^8$ (α,ß) = ((19-ß)/2, (19+ß)/2-α) valid for α = 8, ß = 1

Derives from $e_2^{51}$ x $K_{17}(\alpha_1,ß_1)$ = $(\alpha_1- (ß_1-1)/2, 9- ß_1)$ = $K_{17}(\alpha_2,ß_2)$,

with $e_2^{51}$(α,ß)= [1/2 (α+ß+1), α-1], (8,1) $R_3^8$ (9,2)

This equivalence seems to correspond to $e_2^5$ (α,ß) = (10-ß, 10-α) defined in $A_4$. However, it coincides only in this case and due to

(8,1) $R_2^2$ (2,1) and (9,2) $R_2^3$ (9,8) with (2,1) $R_3^{81}$ (9,8)

The developed equivalences are enough to group classes $C_2$ in classes $C_3$ (Table 9) and establish the transformation flow between classes $C_3$ (Graph 7).

**Equivalences $R_4$ and above**

Similarly to what happened in $A_4$ when increasing the order of an equivalence, new equivalences are produced between specific classes, so the algebraic expressions relating them are not so interesting. The partition in classes corresponding to $R_4$, $R_5$ and $R_6$ are given in Table 9. Graph 7 shows the corresponding transformation trees.

Table 9. Equivalence classes $C_3$, $C_4$, $C_5$ and $C_6$ in $A_5$

| Equivalence | $R_3$ | | | $R_4$ | | | $R_5$ | | | $R_6$ | | |
|---|---|---|---|---|---|---|---|---|---|---|---|---|
| Tree | Classes $C_3$ | Classes $C_2$ included | Image of $C_3$ | Classes $C_4$ | Classes $C_3$ included | Image of $C_4$ | Classes $C_5$ | Classes $C_4$ included | Image of $C_5$ | Classes $C_6$ | Classes $C_5$ included | Image of $C_6$ |
| Tree A | $C_3^1$ | $C_2^1$ | $C_2^2$ | $C_4^4$ | $C_3^1, C_3^5$ | $C_3^4$ | $C_5^3$ | $C_4^4$ | $C_4^3$ | $C_6^3$ | $C_5^3$ | $C_5^2$ |
| | $C_3^2$ | $C_2^3$ | $C_2^1$ | $C_4^1$ | $C_3^2$ | $C_3^4$ | $C_5^4$ | $C_4^1, C_4^5$ | $C_4^4$ | $C_6^4$ | $C_5^4$ | $C_5^3$ |
| | $C_3^3$ | $C_2^4$ | $C_2^3$ | $C_4^2$ | $C_3^3$ | $C_3^2$ | $C_5^1$ | $C_4^2$ | $C_4^1$ | $C_6^1$ | $C_5^1, C_5^5$ | $C_5^4$ |
| | $C_3^4$ | $C_2^2, C_2^5$ | $C_2^4$ | $C_4^3$ | $C_3^4$ | $C_3^3$ | $C_5^2$ | $C_4^3$ | $C_4^2$ | $C_6^2$ | $C_5^2$ | $C_5^1$ |
| | $C_3^5$ | $C_2^6, C_2^7$ $C_2^8, C_2^9$ | $C_2^5$ | | | | | | | | | |
| | $C_3^6$ | $C_2^{10}$ | $C_2^6$ | $C_4^5$ | $C_3^6, C_3^7$ $C_3^8$ | $C_3^5$ | | | | | | |
| | $C_3^7$ | $C_2^{11}$ | $C_2^7$ | | | | | | | | | |
| | $C_3^8$ | $C_2^{12}$ | $C_2^8$ | | | | | | | | | |
| | $C_3^9$ | $C_2^{13}$ | $C_2^{11}$ | $C_4^6$ | $C_3^9$ | $C_3^7$ | $C_5^5$ | $C_4^6$ | $C_4^5$ | | | |
| Tree B | $C_3^{10}$ | $C_2^{14}$ | $C_2^{15}$ | $C_4^7$ | $C_3^{10}$ | $C_3^{11}$ | $C_5^6$ | $C_4^7$ | $C_4^8$ | $C_6^5$ | $C_5^6$ | $C_5^7$ |
| | $C_3^{11}$ | $C_2^{15}$ | $C_2^{14}$ | $C_4^8$ | $C_3^{11}$ | $C_3^{10}$ | $C_5^7$ | $C_4^8$ | $C_4^7$ | $C_6^6$ | $C_5^7$ | $C_5^6$ |
| Tree C | $C_3^{12}$ | $C_2^{18}$ | $C_2^{17}$ | $C_4^{10}$ | $C_3^{12}$ | $C_3^{13}$ | $C_5^{11}$ | $C_4^{10}, C_4^{14}$ | $C_4^{12}$ | $C_6^7$ | $C_5^{11}, C_5^{12}$ | $C_5^{10}$ |
| | $C_3^{13}$ | $C_2^{17}$ | $C_2^{16}$ | $C_4^{12}$ | $C_3^{13}, C_3^{17}$ | $C_3^{15}$ | $C_5^{10}$ | $C_4^{12}, C_4^{13}$ | $C_4^{11}$ | $C_6^8$ | $C_5^{10}$ | $C_5^8$ |
| | $C_3^{14}$ | $C_2^{19}, C_2^{20}$ | $C_2^{18}$ | $C_4^9$ | $C_3^{14}$ | $C_3^{12}$ | $C_5^9$ | $C_4^9$ | $C_4^{10}$ | $C_6^{10}$ | $C_5^9$ | $C_5^{11}$ |
| | $C_3^{15}$ | $C_2^{16}, C_2^{22}$ | $C_2^{20}$ | $C_4^{11}$ | $C_3^{15}, C_3^{16}$ | $C_3^{14}$ | $C_5^8$ | $C_4^{11}$ | $C_4^9$ | $C_6^9$ | $C_5^8$ | $C_5^9$ |
| | $C_3^{16}$ | $C_2^{21}$ | $C_2^{19}$ | | | | | | | | | |
| | $C_3^{17}$ | $C_2^{24}$ | $C_2^{22}$ | | | | | | | | | |
| | $C_3^{18}$ | $C_2^{23}$ | $C_2^{21}$ | $C_4^{13}$ | $C_3^{18}$ | $C_3^{16}$ | | | | | | |
| | $C_3^{19}$ | $C_2^{26}$ | $C_2^{24}$ | $C_4^{14}$ | $C_3^{19}$ | $C_3^{17}$ | | | | | | |
| | $C_3^{20}$ | $C_2^{25}$ | $C_2^{23}$ | $C_4^{15}$ | $C_3^{20}$ | $C_3^{18}$ | $C_5^{12}$ | $C_4^{15}$ | $C_4^{13}$ | | | |

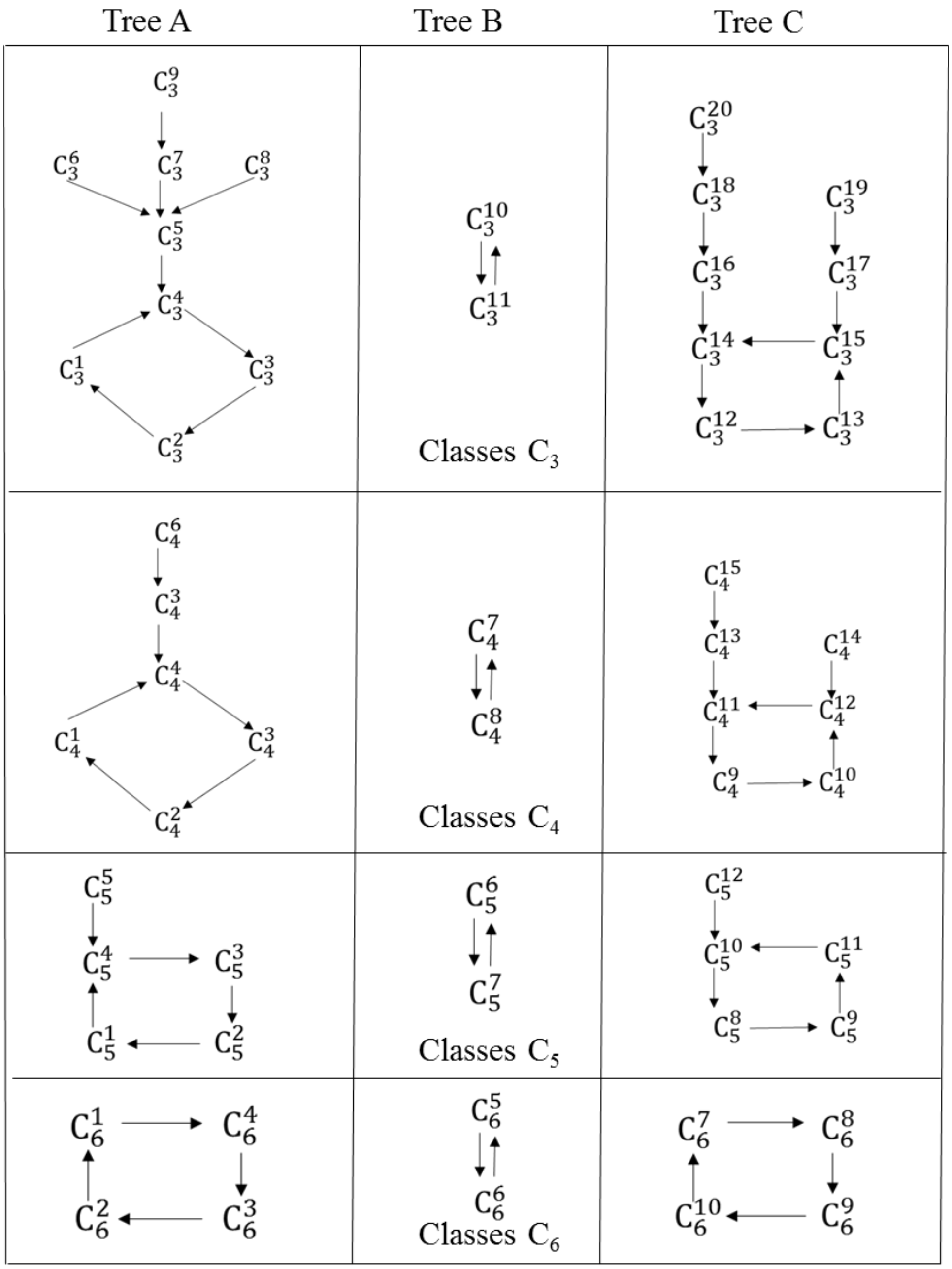

Graph 7. Transformations of classes $C_3$, $C_4$, $C_5$ and $C_6$

## 11. Discussion

In the first part of this paper (Nuez, 2021) we formalized the transformation process by establishing the appropriate functions. The basic functions show that a number's image of w-digits depends on h parameters w=2h+u, u=0, 1, being the reference set $A_w \subset \mathbb{Z}$. In $A_2$ and $A_3$, only one parameter, and two parameters for $A_4$ and $A_5$. All numbers with the same parameters yield the same image. This property has already been established by Prichett et al. (1981).

We developed as well functions $K_i$ that transform such parameters. This contribution is new. Based on these functions, we were able to deduct why transformations lead to 6174 and 495 in $A_4$ and $A_3$, respectively, and to cycles in $A_2$ and $A_5$.

In this second part of the paper we intended to delve into the transformation trees' – or graphs – architecture. We opted for a methodological approach consisting in establishing the algebraic relations that rule the branching of a tree. (secondary branches, main branches…)

To this end, we introduce the notion of equivalences of order r, $R_r$, defined between pairs of numbers whose images match after r transformations [1]. These binary relations satisfy the reflexive, symmetric and transitive features, and they allow to establish partitions of $A_w$ in equivalence classes $C_r$. The main properties of relations $R_r$ and their classes $C_r$, summarized in expressions [2] to [14], are analyzed. All this development is brand new.

A very important part of this are classes $C_1$ which coincide with the sets of numbers that share the same parameters. In $A_2$ and $A_3$, $C_1$ = (α); in $A_4$ and $A_5$, $C_1$ = (α,ß).

How do classes $C_1$ group to create the branches of the tree? Two classes $C_1^i$ and $C_1^j$ whose images match K $(C_1^i)$ = K $(C_1^j)$ will be part of a greater structure $C_2^h = \{C_1^i, C_1^j\}$. In all analyzed sets $A_w$ there is an $R_2^1$ such that for $2 \leq α \leq 9$

$A_2$ and $A_3$: (α) $R_2^1$ (11-α); $A_4$ and $A_5$: (α,0) $R_2^1$ (11-α,0)

For instance, numbers 50595 and 65766 yield the same image after two transformations $K^2$ (50595) = $K^2$ (65766) = 80982 because these numbers' parameters p(50595) = (9,0), p(65766) = (2,0) satisfy $R_2^1$. Note that K (9,0) = K (2,0) = (8,1). This equivalence $α_1 + α_2 = 11$ is valid for any function

$$f(α) = (α\text{-}1 \;\; \overbrace{9 \ldots\ldots 9}^{s} \;\; 10\text{-}α) \quad α \geq 1,\; s \geq 0$$

since the replacement between α-1 and 10-α is compatible with f and it generates a permutation of the same digits. This makes the images of both permutations match and the classes that such permutations create are $(α_1)$ $R_2$ $(α_2 = 11\text{-}α_1)$, $α_1 \geq 2$

This relation of second order, together with the identity (α) $R_2^0$ (α), or its associated functions [14] $e_2^0$ (α) = α, $e_2^1$ (α) = (11-α), built according to the product of functions, has an isomorphic group structure in $\mathbb{Z}_2$.

In $A_4$ and $A_5$ there are new equivalences $R_2$ (group I) $e_2^0$, $e_2^2$, $e_2^3$ *and* $e_2^4$. All of them correspond to transpositions between extreme digits, middle ones or both, in the sequence of the transformed number. These four equivalences have an isomorphic abelian group structure of Klein group.

In $A_4$ there are another four $R_2$ (set II), $e_2^5$, $e_2^6$, $e_2^7$ and $e_2^8$ which do not have a group structure. $e_2^5$ and $e_2^8$ are symmetric, but $e_2^6$ and $e_2^7$ are not. The symmetric failure in the latter is due to a domain-of-existence problem, as a consequence of the restriction α ≥ ß in the images.

However, these four $R_2$ do not exist in $A_5$, except when coinciding with those in group I. This is because there are no two functions $K_i$, $K_j$ yielding the appropriate images.

Functions $K_i$ in $A_4$ and $A_5$ associated to the same permutations – and therefore with identical domains of existence – yield different parameters α' and ß'. The presence of a 9 in the image in $A_5$ by $f_l$ is the cause of such difference. In $A_5$ there are also other equivalences $R_2$ which do not occur in $A_4$.

The new higher structures $C_r$ $(r > 2)$ group different classes $C_{r-1}$ which yield the same image $C_{r-1}$. This enables the analysis of the trees (Tables 1, 5, 7, 8 and 9; Graphs 2, 4, 5, 6 and 7 and those shown in sections 8 and 9).

At an operational level, we studied relations $R_3$ and higher according to [12] and [13]. Hence, we were able to establish algebraic relations between classes $C_1$ which yield identical images $K^r(C_1)$. Particularly, equivalences $R_3$ were developed by using approach C, which implies the existence of functions $K_i$ and $K_j$ in whose domains are classes $(\alpha_1,ß_1)$ and $(\alpha_2,ß_2)$ respectively, as well as an $e_2$ compatible with the image of $K_i(\alpha_1,ß_1)$.

As higher equivalence orders are studied, the domains of existence decrease, so that some $R_r$ can only be applied to a couple of classes $C_1$. This means that their algebraic relation is not that interesting. Such is the case of $R_5^2$ in $A_4$ (2t-19, 28-2t) $R_5$ (2t-21, 30-2t) only valid for t=13 which leads to (7,2) $R_5^2$ (5,4).

A topic of special interest is the relationship of the cycles with the architecture of the tree. There is a cycle of r links if and only if for $(\alpha,ß) = (\alpha_1, ß_1)$ there is a product of operators $K^r = K_q \times \underbrace{\ldots\ldots}_{r} \times K_j \times K_i$ such that $K^r(\alpha_1,ß_1) = (\alpha_1, ß_1)$. In $A_2$ there is only one cycle of 5 links; in $A_5$ there is one 2-link cycle and two 4-link cycles.

The existence of cycles affects the structure of the transformation trees. In $A_2$ there is simple branching, of a single class, on the structure of the cycle. In $A_5$ there are 3 trees: tree A contains a cycle formed by classes (7,5), (6,3), (7,2) and (8,4), and it groups 30 classes (α,ß); tree B, with a cycle formed by (6,0) and (5,4), groups 3 classes $C_1$; tree C, with the cycle (7,6), (6,4), (6,2) and (8,3), groups 21 classes $C_1$.

Is there any relationship between trees and functions $K_i$? All functions $K_i$ operate on classes (α,ß) of the two main trees A and C in $A_5$, except for functions $K_{13}$, $K_{21}$, $K_{25}$ and $K_{26}$. $K_{25}$ and $K_{26}$ only operate on classes (α,0) and they are found in trees A and B which have these classes; in tree C, there are not classes (α,0). $K_{13}$ only operates on (5,4) and (5,5), and it is found in tree C. Eight out of the 13 functions $K_i$ operate on tree B, but this tree only has 3 classes. Therefore, such relationship does not seem to occur.

What about relationships between trees and combinations of functions $K_i$? It is interesting to note that the same combination of 4 functions exists between the functions operating within the main cycles A and C

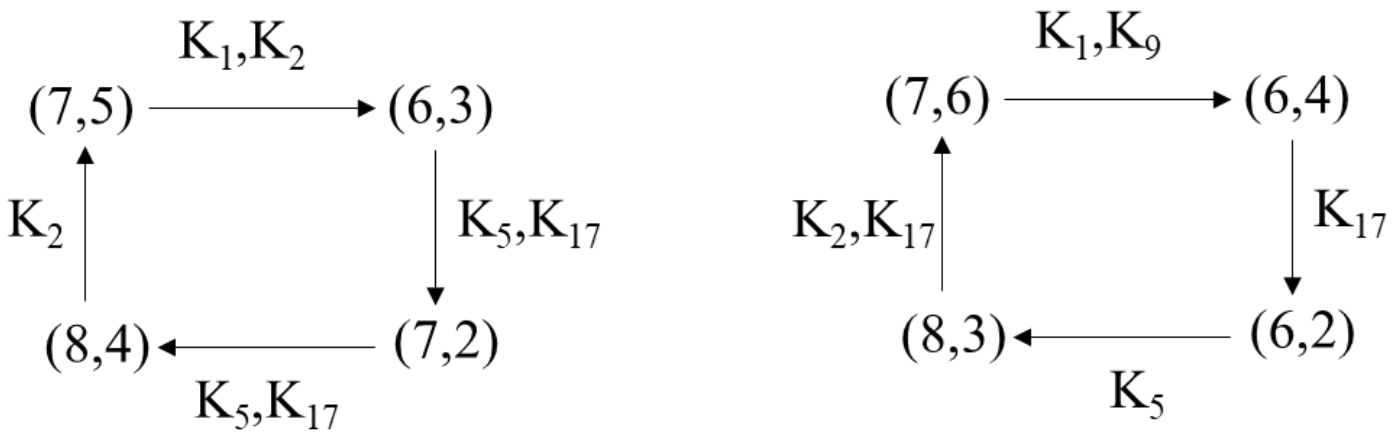

The shared sequence of operators is $K_5$, $K_{17}$, $K_1$, $K_2$ so that for

$K^4 = K_5 \times K_{17} \times K_1 \times K_2$, $K^4\,(8,4) = (8,4)$, $K^4\,(8,3) = (8,3)$

$K^4 = K_2 \times K_5 \times K_{17} \times K_1$, $K^4\,(7,5) = (7,5)$, $K^4\,(7,6) = (7,6)$

and analogously in the other two cases of functions $K_i$ rotation.

Moreover, the main branches in tree A are articulated through $(8,6) \xrightarrow{K_1} (7,5)$, while in tree C, they are so through $(8,7) \xrightarrow{K_1,K_9} (7,6)$.

It is evident that the same combinations of operators $K_i$, when acting on different classes, generate separate structures. It is the combination of operators and domains of existence –classes (α,ß)– that determines the structure. Such combination is a consequence of basic functions $f_i$ changing with each $A_w$.

There are as many new equivalence classes –[5] and [5a], from the last level– as there are nodes in a cycle. In $A_5$, tree A has four ($C_6^1$, $C_6^2$, $C_6^3$ and $C_6^4$); tree B has two ($C_6^5$ and $C_6^6$) and tree C has four ($C_6^7$, $C_6^8$, $C_6^9$ and $C_6^{10}$). These classes contain all those classes of a lower level that are as many steps away from a cycle's nodal class as there are links in such cycle, or a multiple of that number. Thus, $C_6^4$ is made up of 4 classes $C_3$ ($C_3^2$, $C_3^6$, $C_3^7$ and $C_3^8$), from which the nodal class is $C_3^2$. Naturally, being $C_3^2$ situated in the cycle, it requires 4 steps or transformations to go back to $C_3^2$. The other classes ($C_3^6$, $C_3^7$ and $C_3^8$) require each 4 steps to become $C_3^2$. The same applies for the rest of the classes $C_5$, $C_4$, $C_2$ or $C_1$. For example, let $C_6^4$ contain classes (4,0), (7,0), (7,5), (1,1), (9,1), (9,9), (2,0), (9,0) and (6,5). The cycle's nodal class is (7,5). Any other class needs 4 transformations for its image to be (7,5).

In tree B with a 2-link cycle, $C_6^5 = C_2^{14}$ contains classes (5,0) and (6,0), while $C_6^6 = C_2^{15}$ only class (5,4). The cycle's nodal classes are (6,0) and (5,4). (5,0) requires two steps to become (6,0), just like (6,0) itself or (5,4).

In the case of $A_2$ the cycle has 5 links, which is the number of steps required by each of the classes (α) embedding classes $C_2$ to go back to the cyclic node of their class.

If the branches of the tree are far away enough from the cycle, there may be classes which are more steps away from entering the cycle than the number of links in the cycle itself. Such is the case of (1,0), which along with (8,4) and (8,6) make up class $C_6^1$ in $A_5$'s tree A. Their distance to the nodal class (8,4) is 8, which is a multiple of 4. **The length of the cycle** –its number of links– **is the pattern of distances which determines how classes are grouped** around the cycle's nodes. This feature is a direct consequence of definitions [1] and [7].

The transformation trees with a single balance constant, as in $A_4$ and $A_3$, with balance classes (6,2) and (5) respectively, are actualy structures with cycles of a single link (6,2), (5). Consequently, there can only be one higher equivalence class – $C_7^1$ in $A_4$ and $C_6^1$ in A3. Furthermore, such equivalence class embeds all classes of a lower order, since any number of steps is a multiple of 1, which is the number of links in the cycle. As previously indicated when studying $A_3$'s tree, in each increase to a higher equivalence all the previous equivalences stay [8], and the class containing the balance class merges with the adjoining one (Graph 4 and Table 7 for $A_4$).

Hence, it turns out that the methodological approach C (Graph 3), based on a hierarchical and inclusive equivalence system, converges with approach B of inclusive layers.

**Acknowledgement**

To the UPV professor, Dr. María José Díez, for her selfless collaboration in correcting the manuscript.

ANNEX

Some tables from the first part of this paper can be found in this annex

Table 2. Functions $K_i$ (α,ß) = (α',ß') in $A_4$ with ß > 0

| Type P= O(n') | | $K_i$ (α,ß) | (α', ß') | Existence conditions |
|---|---|---|---|---|
| 1a | $P_1$ (1 2 3 4) | $K_1$ | (2α−10, 2ß-10) | α≥ß+1, ß≥5 |
| | $P_2$ (1 3 2 4) | $K_2$ | (2α−10, 10-2ß) | α≥6, ß≤5, α+ß≥11 |
| 1b | $P_5$ (3 1 4 2) | $K_5$ | (10-2ß, 2α−10) | α≥5, α+ß≤9 |
| | $P_6$ (3 4 1 2) | $K_6$ | (10-2ß, 10−2α) | α≤5, α≥ß+1 |
| 2a | $P_9$ (1 2 4 3) | $K_9$ | [(α+ß)-9, (α+ß)]-11] | α≤ ß+1, α+ß≥11 |
| | $P_{10}$ (1 4 2 3) | $K_{10}$ | [(α+ß)-9, 11-(α+ß)] | α≥5, ß≥5, α+ß≤11 |
| 2b | $P_{13}$ (4 1 3 2) | $K_{13}$ | [11-(α+ß), (α+ß)-9] | α≤5, ß≤5, α+ß≥9 |
| | $P_{14}$ (4 3 1 2) | $K_{14}$ | [11-(α+ß), 9-(α+ß)] | α≤ß+1, α+ß≤9 |
| 3a | $P_{17}$ (1 3 4 2) | $K_{17}$ | [(α-ß)+1, (α-ß)-1] | α≥ß+1, 9≤ α+ß≤11 |
| | $P_{18}$ (1 4 3 2) | $K_{18}$ | [(α-ß)+1, 1-(α-ß)] | α≥5, ß≤5, α≤ß+1 |
| 3b | $P_{21}$ (4 1 2 3) | $K_{21}$ | [1-(α-ß), (α-ß)+1] | α= ß=5 |

Table 3. Functions $K_i$ (α,0) = (α',ß') in $A_4$

| P = O(n') | $K_i$ | (α',ß') | Existence conditions |
|---|---|---|---|
| $P_{25}$(1 2) | $K_{25}$ (α,0) | (α-1, 10-α) | 6≤α≤9 |
| $P_{26}$ (2 1) | $K_{26}$ (α,0) | (10-α, α-1) | 0<α≤5 |

Table 4. Functions $K_i$ in $A_5$

| Type | P=O(n') | Existence conditions | $K_i(\alpha,\beta) = (\alpha',\beta')$ |
|---|---|---|---|
| | $\beta>0$ $P_{1234} = (9\ \ \alpha\ \ \beta-1\ \ 9-\beta\ \ 10-\alpha)$, $\alpha>\beta$ | | |
| 1a | $P_1$ (1 2 3 4) | $\alpha\geq\beta+1$, $5\leq\beta\leq 8$ | $K_1(\alpha,\beta) = (\alpha-1,\alpha+\beta-9)$ |
| | $P_2$ (1 3 2 4) | $\alpha\geq 6$, $2\leq\beta\leq 5$, $\alpha+\beta\geq 11$ | $K_2(\alpha,\beta) = [\alpha-1,(\alpha-\beta)+1]$ |
| 2a | $P_5$ (3 1 4 2) | $\alpha\geq 5$, $\alpha+\beta\leq 9$ | $K_5(\alpha,\beta) = [10-\beta,(\alpha-\beta)-1]$ |
| | $P_6$ (3 4 1 2) | $\alpha\leq 5$, $\alpha\geq\beta+1$ | $K_6(\alpha,\beta) = [10-\beta,9-(\alpha+\beta)]$ |
| 2a | $P_9$ (1 2 4 3) | $\alpha\leq\beta+1$, $\alpha+\beta\geq 11$ | $K_9(\alpha,\beta) = (\beta,2\alpha-10)$ |
| | $P_{10}$ (1 4 2 3) | $\alpha\geq 5$, $\beta\geq 5$, $\alpha+\beta\leq 11$ | $K_{10}(\alpha,\beta) = [\beta,(\alpha-\beta)+1]$ |
| 2b | $P_{13}$ (4 1 3 2) | $\alpha\leq 5$, $\beta\leq 5$, $\alpha+\beta\geq 9$ | $K_{13}(\alpha,\beta) = [10-\beta,1-(\alpha-\beta)]$ |
| | $P_{14}$ (4 3 1 2) | $\alpha\leq\beta+1$, $\alpha+\beta\leq 9$ | $K_{14}(\alpha,\beta) = (10-\beta,10-2\alpha)$ |
| 3a | $P_{17}$ (1 3 4 2) | $\alpha\geq\beta+1$, $9\leq\alpha+\beta\leq 11$ | $K_{17}(\alpha,\beta) = (10-\beta,2\alpha-10)$ |
| | $P_{18}$ (1 4 3 2) | $\alpha\geq 5$, $\beta\leq 5$, $\alpha\leq\beta+1$ | $K_{18}(\alpha,\beta) = [10-\beta,(\alpha+\beta)-9]$ |
| 3b | $P_{21}$ (4 1 2 3) | $\alpha=\beta=5$ | $K_{21}(\alpha,\beta) = [\beta,11-(\alpha+\beta)]$ |
| | $\beta=0$ $P_{12} = (9\ \ 9\ \ 9\ \ \alpha-1\ \ 10-\alpha)$ | | |
| 4 | $P_{25}$ (1 2) | $6<\alpha\leq 9$ | $K_{25}(\alpha,0) = (\alpha-1,10-\alpha)$ |
| | $P_{26}$ (2 1) | $0<\alpha\leq 5$ | $K_{26}(\alpha,0) = (10-\alpha, \alpha-1)$ |